\documentclass[11pt,a4paper,reqno]{amsart}

\usepackage{amscd,amsmath}
\usepackage{amssymb,amsmath,latexsym,color,enumerate,dsfont}
\usepackage[all]{xypic}       
\usepackage{mathabx}
\usepackage{mathrsfs,graphicx}
\usepackage{stmaryrd}

\usepackage{amsthm}
\swapnumbers

\usepackage{tikz}
\usetikzlibrary{matrix}

\def\Pt{{\Sigma}}
\def\O{{\mathcal O}}

\def\LR{{\mathfrak{L}(\R)}}

\def\intervalR{\mathbb{IR}}
\def\N{{\mathds{N}}}
\def\R{{\mathbb{R}}}
\def\Q{{\mathbb{Q}}}

\def\LR{{\mathfrak{L}(\R)}}

\def\N{{\mathbb{N}}}

\def\R{{\mathbb{R}}}
\def\Q{{\mathbb{Q}}}

\def\Frm{{\mathsf{Frm}}}
\def\Top{{\mathsf{Top}}}

\newcommand{\newuparrow}{{{\rlap{$\ $}\hbox{$\uparrow$}}}}%
\newcommand{\twoheaddownarrow}{{\rlap{\rlap{$\ $}\raise .25ex\hbox{$\downarrow$}}\raise-.25ex\hbox{$\downarrow$}}}
\newcommand{\twoheaduparrow}{{\rlap{\rlap{$\ $}\raise .25ex\hbox{$\uparrow$}}\raise-.25ex\hbox{$\uparrow$}}}

 \newcommand{\tbigcup}{\mathop{\textstyle \bigcup }}
 \newcommand{\tbigvee
}{\mathop{\textstyle \bigvee }}
 \newcommand{\tbigwedge
}{\mathop{\textstyle \bigwedge }}

 \newcommand{\tbigcap
}{\mathop{\textstyle \bigcap }}

 \newcommand{\gidoia}{\textsf{\scriptsize ---}}

\newcommand{\vu}{\mathop{\vee\hskip-6.7pt \wedge}}

\DeclareMathAccent{\widetriangle}{\mathord}{largesymbols}{"E6}

\numberwithin{equation}{section}

\theoremstyle{plain}
\newtheorem{theorem}{Theorem}[subsection]
\newtheorem{lemma}[theorem]{Lemma}
\newtheorem{proposition}[theorem]{Proposition}
\newtheorem{corollary}[theorem]{Corollary}

\theoremstyle{definition}

\newtheorem{remark}[theorem]{Remark}
\newtheorem{remarks}[theorem]{Remarks}

\newtheorem{example}[theorem]{Example}
\newtheorem{examples}[theorem]{Examples}

\newtheorem*{factsn}{Facts}
\newtheorem*{remarkn}{Remark}

\newtheorem*{examplesn}{Examples}

\newcommand\coS{\reflectbox{$\mathsf{S}$}}

\begin{document}

\title[]{A lattice-theoretic approach\\
to arbitrary real functions on frames}

\author[I. Mozo Carollo]{Imanol~Mozo~Carollo}

\address{Department of Mathematics\\
University of the Basque Country UPV/EHU\\Apdo. 644, 48080 Bilbao\\Spain}
\address{CECAT - Department of Mathematics \& Computer Science\\ Chapman University\\ One University Drive, Orange, CA  92866, USA}
\email{imanol.mozo@ehu.eus}

\date{\today}

\thanks{
I. Mozo Carollo acknowledges financial support from the Ministry of Economy and Competitiveness of Spain under grant MTM2015-63608-P (MINECO/FEDER) as well as from the Basque Government under grant IT974-16 and Postdoctoral Fellowship POS-2015-1-0103.\\
The Version of Record of this manuscript has been published online on 07 Nov 2017 and is available in Quaestiones Mathematicae, DOI: 10.2989/16073606.2017.1380725.}

\subjclass[2010]{Primary: 06D22; Secondary:
26A15, 54C30, 54D15.}

\keywords{frame, locale, frame of reals,
 continuous real function,  order complete, Dedekind-MacNeille completion, semicontinuous real function, partial real function, Hausdorff continuous real function.}

\begin{abstract} 
In this paper we model discontinuous extended real functions in pointfree topology following a lattice-theoretic approach, in such a way that, if $L$ is a subfit frame, arbitrary extended real functions on $L$ are the elements of the Dedekind-MacNeille completion of the poset of all extended semicontinuous functions on $L$. This approach mimicks the situation one has with a $T_1$-space $X$, where the lattice $\overline{\mathrm{F}}(X)$ of arbitrary extended real functions on $X$ is the smallest complete lattice containing both extended upper and lower semicontinuous functions on $X$. Then, we identify real-valued functions by lattice-theoretic means. By construction,  we obtain definitions of discontinuous functions that are conservative for $T_1$-spaces. We also analyze semicontinuity and introduce definitions which are conservative for $T_D$-spaces.
\end{abstract}

\maketitle

\section*{Introduction}

In pointfree topology frames are viewed as generalized spaces and frame homomorphism as continuous functions between them.  Indeed, if $X$ is a sober space, continuous functions $Y\to X$ are in one-to-one correspondence with frame homomorphism $\mathcal{O}X\to\mathcal{O}Y$. Motivated by the fact that arbitrary real functions on a space $X$ are continuous functions when one replaces $\mathcal{O}X$ by the discrete topology, which can be thought of as the system of subspaces of $X$, Guti\'errez Garc\'\i a, Kubiak and Picado defined arbitrary real functions on a frame $L$ as  frame homomorphisms from the frame of reals $\mathfrak{L}(\R)$ (the pointfree counterpart of the real line) to the dual $\coS(L)=\mathsf{S}(L)^\mathrm{op}$ of the system of all sublocales of $L$ \cite{GKP2009}. Certainly this was shown to be successful definition that provides a very nice setting for the treatment of semicontinuity. It made possible to extend many classical definitions, constructions and results to the pointfree setting \cite{GKP09a,GKP09c,GK11,GMP14}. Recently, in \cite{PP16} Picado and Pultr have replaced $\coS(L)$  by a smaller class of sublocales (see \cite{PPT16}), mending some disadvantages that the former definition had. For example, this new approach yields a ``lattice of arbitrary real functions'' which is Dedekind complete in \cite{PP16}, while, in general, it was not so in \cite{GKP2009}.

Beyond the differences between \cite{GKP2009} and \cite{PP16}, the idea behind both definitions is the same: replacing $L$ by some discretization of it which will play the role of the discrete topology in the classical theory. The main goal of this paper is to approach this problem from a different point of view. The motivation comes from the fact the lattice $\overline{\mathrm{F}}(X)$ of arbitrary extended real functions on a $T_1$-space $X$ is the smallest complete lattice in which the poset of lower and upper semicontinuous extended real functions embeds, more rigorously, $\overline{\mathrm{F}}(X)$ is the Dedekind-MacNeille completion of the poset $\overline{\mathrm{LSC}}(X)\cup\overline{\mathrm{USC}}(X)$ of all semicontinuous extended real functions on $X$. For this purpose, we will need a proper description of extended semicontinuous functions on frames and a convenient representation in terms of certain kind of functions of the Dedekind-MacNeille completion of poset they form. We should note that, even though in \cite{GKP2009} and \cite{PP16} the authors focused on real-valued functions, their approach can be immediately applied to arbitrary frames, while ours is intrinsically linked to the order structure of the reals. However, due to the importance of the real line and real-valued functions in topology, we believe that our approach is worth to be investigated.

There are various equivalent ways of introducing the \emph{frame of reals} $\LR$ \cite{BB97}.
Here it will be useful to adopt the description used in \cite{BGP14} given by generators $(r,\gidoia)$ and $(\gidoia,s)$, $r,s\in\Q$,  and defining relations:
 \begin{enumerate}[\rm(r1)]
 \item $(r,\gidoia)\wedge (\gidoia,s)=0$ whenever $r\ge s$,
 \item $(r,\gidoia)\vee (\gidoia,s)=1$ whenever $r<s$,
 \item $(r,\gidoia)= \tbigvee_{s>r}(s,\gidoia)$, for every $r\in{{\mathbb{Q}}}$,
 \item $(\gidoia,s)= \tbigvee_{r<s}(\gidoia,r)$, for every $r\in {{\mathbb{Q}}}$,
 \item $ \tbigvee_{r\in \mathbb{Q}}(r,\gidoia)=1$,
 \item $\tbigvee_{s\in \mathbb{Q}}(\gidoia,s)=1$.
 \end{enumerate}
 By dropping relations $(\mathrm{r5 })$ and $(\mathrm{r6})$ in the description of ${{\mathfrak{L}({{\mathbb{R}}})}}$ above, we have the corresponding \emph{frame of extended reals} ${\mathfrak{L}\big(\overline{{{\mathbb{R}}}}\big)}$ \cite{BGP14}. For any frame $L$, a \emph{continuous real function} \cite{BB97} (resp. \emph{extended continuous real function} \cite{BGP14}) on a frame $L$ is a frame homomorphism $f\colon {\mathfrak{L}(\mathbb{R})}\rightarrow L$ (resp. $f\colon {\mathfrak{L}\big(\overline{{{\mathbb{R}}}}\big)}\rightarrow L$). We denote by $\mathrm{C}(L)$ (resp. $\overline{\mathrm{C}}(L)$) the collection of all (resp. extended) continuous real functions on $L$. The correspondences $L\mapsto \mathrm{C}(L)$ and $L\mapsto \overline{\mathrm{C}}(L)$ are functorial in the obvious way.

Extended lower semicontinuous functions $X\to\overline{\R}$ are in one-to-one correspondence with frame homomorphisms $\mathfrak{L}_u(\overline{\R})\to\mathcal{O}X$, where $\mathfrak{L}_u(\overline{\R})$ is a convenient subframe of $\mathfrak{L}(\overline{\R})$. One can represent upper semicontinuous functions similarly by replacing $\mathfrak{L}_u(\overline{\R})$ by another convenient subframe $\mathfrak{L}_l(\overline{\R})$ of $\mathfrak{L}(\overline{\R})$. Replacing $\mathcal{O}X$ by an arbitrary frame $L$ provides a natural extension of the notion of semicontinuous functions. As shown in \cite{GKP2009}, we can represent both extended semicontinuous functions on a frame $L$ as frame homomorphisms $\mathfrak{L}(\overline{\R})\to\coS(L)$. We define in this fashion 
 the class $\overline{\mathrm{LSC}}(L)$ resp. $\overline{\mathrm{USC}}(L)$ of lower resp. upper semicontinuous functions on $L$ (Section 2).

Their non-extended versions are not that simple to define. Semicontinuous functions were first introduced by Li and Wang in \cite{LW97}, but their definition did not represent the classical notion. It was mended in \cite{GP07} and \cite{GKP2009} in a way that made possible to formulate Kat\u etov-Tong insertion theorem in the pointfree setting. However, this definition fails  to represent all semicontinuous functions from the classical theory. This small inaccuracy has already been pointed out and corrected in \cite{PP16}. In Section 2, we provide a definition in the setting that fits in our approach to arbitrary real functions. In addition, we show that our definition extends faithfully the classical notion for all $T_D$-spaces, while the one in \cite{PP16} was shown to do so only for $T_1$-spaces.

In Section 3, we derive a description of the Dedekind-MacNeille completion of the lattice $\overline{\mathrm{C}}(L)$ from the Dedekind completion of the lattice $\mathrm{C}^\ast(L)$ of bounded continuous functions on $L$ in terms of continuous partial real functions presented in \cite{MGP14}. The frame $\mathfrak{L}(\overline{\mathbb{IR}})$ of \emph{extended partial real numbers} is the frame that we obtain after removing (r2) from the definition of $\mathfrak{L}(\overline{\R})$ and $\overline{\mathrm{IC}}(L)$ the class of all frame homomorphisms $\mathfrak{L}(\overline{\mathbb{IR}})\to L$, namely, \emph{continuous extended partial real functions} on $L$, in which $\overline{\mathrm{C}}(L)$ embeds canonically.  Given a completely regular frame $L$, we describe the Dedekind-MacNeille completion of $\overline{\mathrm{C}}(L)$ inside $\overline{\mathrm{IC}}(L)$ by
\[
\overline{\mathrm{H}}(L)=\{f\in\overline{\mathrm{IC}}(L)\mid f(r,\gidoia)^\ast\leq f(\gidoia,s)\text{ and }f(\gidoia,s)^\ast\leq f(r,\gidoia)\text{ if } r<s\}.
\]

In Section 4, we show that given a subfit frame $L$, the poset $\overline{\mathrm{LSC}}(L)\cup\overline{\mathrm{USC}}(L)$ is join- and meet-dense in $\overline{\mathrm{C}}(\coS(L))$ and consequently their Dedekind-MacNeille completions coincide. This is the motivation to consider the latttice $\overline{\mathrm{F}}(L)=\overline{\mathrm{H}}(\coS(L))$, that is, the Dedekind-MacNeille completion of $\overline{\mathrm{C}}(\coS(L))$, as the pointfree counterpart of the $\overline{\mathrm{F}}(X)$. Accordingly, an \emph{arbitrary extended real function} on a frame $L$ will be a frame homomorphism $f\colon\mathfrak{L}(\overline{\mathbb{IR}})\to\coS(L)$ such that
\[
f(r,\gidoia)^\ast\leq f(\gidoia,s)\text{ and }f(\gidoia,s)^\ast\leq f(r,\gidoia)\quad\text{whenever}\quad r<s.
\]
This condition holds for all $f\in\overline{\mathrm{C}}(\coS(L))$. In comparison with \cite{GKP2009}, an arbitrary extended real function has to satisfy a slightly weaker condition than preserving (r2).

In order to identify real-valued functions among arbitrary extended real functions, one might be tempted to simply replace $\mathfrak{L}(\overline{\mathbb{IR}})$ by its non-extended version, the frame $\mathfrak{L}(\mathbb{IR})$ of partial real numbers \cite{MGP14}, and say that an $f\in\overline{\mathrm{F}}(L)$ is real-valued if it preserves (r5) and (r6). However, this is too restrictive and it does not yield the proper notion. In fact, not even all semicontinuous would be represented with such a definition. In Section 5, we define the class $\mathrm{F}(L)$ of \emph{arbitrary real functions} of $L$ as a subset of $\overline{\mathrm{F}}(L)$ by lattice-theoretic means and characterize them as follows: an $f\in\overline{\mathrm{F}}(L)$ is \emph{real-valued} if and only if
\[
\left(\tbigvee_{r\in\Q}f(r,\gidoia)\right)^\ast=0=\left(\tbigvee_{s\in\Q}f(\gidoia,s)\right)^\ast.
\]
Note that this condition is weaker that requiring that $f$ preserves (r5) and (r6). Indeed, instead of having  $\tbigvee_{r\in\Q}f(r,\gidoia)=1=\tbigvee_{s\in\Q}f(\gidoia,s)$, we only have that those joins are dense elements of the frame $\coS(L)$. This approach provides a faithful extension of the classical notions: given a $T_1$-space $X$, one has
\[
\overline{\mathrm{F}}(X)\simeq\overline{\mathrm{F}}(\mathcal{O}X)\quad\text{and}\quad\mathrm{F}(X)\simeq\mathrm{F}(\mathcal{O}X).
\]
We also show that $\mathrm{F}(L)$ is isomorphic to $\mathrm{C}(\mathfrak{B}(\coS(L)))$, where $\mathfrak{B}(\coS(L))$ denotes the Booleanization of $\coS(L)$ \cite{BP96}. This has two important consequences. First, our definition is equivalent to the one presented in \cite{PP16} for subfit frames. And second,  $\mathrm{F}(L)$ can be equipped with a lattice-ordered ring structure that making possible to consider a theory of rings of real functions.

Finally, in Section 6, we study when our definitions coincides with the one presented in \cite{GKP2009}, that is, when one has $\overline{\mathrm{F}}(L)=\overline{\mathrm{C}}(\coS(L))$ and $\mathrm{F}(L)=\mathrm{C}(\coS(L))$.

\section{Preliminaries}

 \subsection{Dedekind completion}
 For any subset $A$ of a partially ordered set $(P,\le)$, we will denote by $
 {\tbigvee^P} A$ (resp. $
 {\tbigwedge^P} A$) the supremum (resp. infimum) of $A$ in $P$ in case it exists (we shall omit the superscript if it is clear from the context).

 A \emph{Dedekind-MacNeille completion} (also called \emph{completion by cuts},  \emph{normal completion} or just \emph{MacNeille completion}) of a poset $P$ is a join- and meet-dense embedding $\varphi\colon P\to M(P)$ in a complete lattice (as usual, a map
$\varphi\colon P \to M(P)$ is said to be \emph{join-dense} if and only if so is its image in $M(P)$; that is, each element of $M(P)$
is a join of elements from $\varphi[P]$; \emph{meet-density} is defined dually). The Dedekind-MacNeille completion is the only complete lattice in which the given poset is join- and meet-dense.

Sometimes a weaker kind of completion is more userful: a poset $(P,\le)$ is \emph{Dedekind} (\emph{order}) \emph{complete}\index{Dedekind (order) complete}\index{complete! Dedekind (order) $\sim$} (or \emph{conditionally complete}\index{conditionally complete}\index{complete!conditionally $\sim$}) if every non-void subset $A$ of $P$ which is bounded from above has a supremum in $P$ (and then, in particular, every non-void subset $B$ of $P$ which is bounded from below will have a infimum in $P$). Of course, being complete is equivalent to Dedekind complete plus the existence of top and bottom elements.
A \emph{Dedekind completion} (or \emph{conditional completion}) of $P$ is a join- and meet-dense embedding $\varphi\colon P\to D(P)$ in a Dedekind complete poset $D(P)$. The Dedekind completion is slightly smaller than the MacNeille completion: it can be obtained from $M(P)$, in case $P$ is directed, just by removing its top and bottom elements.

For more information on universal properties of the Dedekind-MacNeille and the Dedekind completion see \cite[1.3]{SS10}.

\subsection{Frames}

A \emph{frame}\index{frame} (or \emph{locale}\index{locale}) $L$ is a complete lattice such that
 \begin{equation}
 a\wedge \tbigvee B=\tbigvee \{a\wedge b\mid b\in B\}
 \end{equation}
 for all $a\in L$ and $B\subseteq L$; equivalently, it is a complete Heyting algebra with Heyting operation $\rightarrow$ satisfying the standard equivalence $a\wedge b\le c$ if and only if $a\le b\rightarrow c$. The \emph{pseudocomplement}\index{pseudocomplement} of an $a\in L$ is the element
 \[
 a^{\ast}=a\rightarrow 0=\tbigvee \{b\in L\mid a\wedge b=0\}.
 \]
An element $a\in L$ is \emph{complemented} if $a\vee a^\ast=1$. For completemented elements $a\in L$ the dual distributivity law also holds:
\[
a\vee\tbigwedge B=\tbigwedge\{a\vee b\mid b\in B\}
\]
for all $B\subseteq L$. An element $a\in L$ is \emph{dense}\index{regular!element} if $a^\ast=0$ (equivalently, if $a^{\ast\ast}=1$). A \emph{frame homomorphism}\index{frame homomorphism} is a map $h\colon L\to M$ between frames which preserves finitary meets (including the top element 1) and arbitrary joins (including the bottom element 0). Then $\Frm$ is the corresponding category of frames and their homomorphisms.

 The most typical example of a frame is the lattice $\O X$ of open subsets of a topological space $X$. The correspondence $X\mapsto \O X$ is clearly functorial, and consequently we have a contravariant functor $\O \colon \Top\to\Frm$ where $\Top$ denotes the category of topological spaces and continuous maps. There is also a functor in the opposite direction, the \emph{spectrum functor}\index{spectrum functor}\index{functor!spectrum $\sim$}\index{functor!$\Pt$} $\Sigma\colon\Frm\to\Top$ which assigns to each frame $L$ its spectrum $\Sigma L$, the space of all homomorphisms $\xi\colon L\to \{0,1\}$ with open sets $\Sigma_a=\{\xi\in\Sigma L\mid \xi(a)=1\}$ for any $a\in L$, and to each frame homomorphism $h\colon L\to M$ the continuous map $\Sigma h\colon\Sigma M\to\Sigma L$ such that $\Sigma h(\xi)=\xi h$. The spectrum functor is right adjoint to $\mathcal O$, with adjunction maps $\eta_L\colon L\to \O \Sigma L, \;\eta_L(a)=\Sigma_a$ and $\epsilon_X\colon X\to \Sigma
 \O X, \; \epsilon_X(x)=\hat{x}, \; \hat{x}(U)=1\mbox{ if and only if }x\in U$ (the former is the \emph{spatial reflection} of the frame $L$). A frame is said to be \emph{spatial}\index{spatial frame}\index{frame!spatial $\sim$} if it is isomorphic to the frame of open sets of a space.

For general notions and results concerning frames we refer to Johnstone \cite{PJ82} or the recent
Picado-Pultr \cite{PP12}. The particular notions we will need are the following:
a frame $L$ is:

 \begin{enumerate}[-]
%

\item \emph{subfit} if for each $a,b\in L$ such that $a\not\leq b$ there exists $c\in L$ such that $a\vee c=1\neq b\vee c$;

 \item \emph{completely regular}\index{completely regular frame}\index{frame!completely regular $\sim$} if $a = \tbigvee\{b\in L\mid b{\prec\hskip-2pt \prec\hskip2pt}a\}$ for each $a\in L$, where $b{\prec\hskip-2pt \prec\hskip2pt}a$ ($b$ is \emph{completely below}\index{completely below relation}\index{relation!completely below $\sim$} $a$) means that there is $\{c_{r}\mid r\in \mathbb{Q}\cap[0,1]\}\subseteq L$ such that $b\le c_{0}$, $c_{1}\le a$ and $c_r\prec c_s$ (i.e. $c_{r}^\ast\vee c_{s}=1$) whenever $r<s$;

\item \emph{extremally disconnected} if $a^\ast\vee a^{\ast\ast}=1$ for every $a\in L$.

 
 
  
 \end{enumerate}

\subsection{Sublocales}
 A \emph{sublocale set} (briefly, a
\emph{sublocale}) $S$ of a locale $L$ is a subset $S \subseteq L$ such
that
\begin{itemize}
\item[(S1)] for every $A \subseteq S$, $\bigwedge A$ is in $S$, and
\item[(S2)] for
every $s \in S$ and every $x \in L$, $x \to s$ is in $S$.
\end{itemize}

The system of all sublocales constitutes a coframe with the order given by inclusion, meet coinciding with the intersection and the join given by $\tbigvee S_i=\{\tbigwedge M\mid M\subseteq \tbigcup S_i\};$ the top is $L$ and the bottom is the set $\{1\}$.

\subsubsection{Closed and open sublocales} For any $a\in L$, the sets $\mathfrak{c}(a)=\newuparrow a$ and $\mathfrak{o}(a)=\{a\rightarrow b\mid b\in L\}$ are the \emph{closed} and \emph{open}\index{open sublocale $\sim$} \index{sublocale!open} sublocales of $L$, respectively. The map $a\mapsto \mathfrak{c}(a)$ is a coframe embedding $L^{\mathrm{op}} \hookrightarrow {\mathsf{S}(L)}$ providing an isomorphism $\mathfrak{c}$ between $L^{\mathrm{op}}$ and the subcoframe $\mathfrak{c}(L)$ of ${\mathsf{S}(L)}$ consisting of all closed sublocales. On the other hand, denoting by $\mathfrak{o}(L)$ the subcoframe of ${\mathsf{S}(L)}$ generated by all $\mathfrak{o}(a)$, the correspondence $a\mapsto \mathfrak{o}(a)$ establishes  poset embedding $L\rightarrow\mathfrak{o}(L)$. The following holds:

\begin{factsn} The following facts hold in $\mathsf{S}(L)$.
\begin{enumerate}[(1)]
\item  $\mathfrak{o}(a)$ and $\mathfrak{c}(a)$ are complements of each other.
\item $\mathfrak{c}(0)=1$, $\mathfrak{c}(1)=0$, $\mathfrak{c}(a)\vee\mathfrak{c}(b)=\mathfrak{c}(a\wedge b)$ and $\tbigwedge_{i\in I}\mathfrak{c}(a_i)=\mathfrak{c}\left(\tbigvee_{i\in I}a_i\right)$.
\item $\mathfrak{o}(0)=0$, $\mathfrak{o}(1)=1$, $\mathfrak{o}(a)\wedge\mathfrak{o}(b)=\mathfrak{o}(a\wedge b)$ and $\tbigvee_{i\in I}\mathfrak{o}(a_i)=\mathfrak{o}\left(\tbigvee_{i\in I}a_i\right)$.
\item Each sublocale $S\in\mathsf{S}(L)$ can be represented as $S=\tbigwedge_i(\mathfrak{o}(a_i)\vee\mathfrak{c}(b_i))$.
\end{enumerate}
\end{factsn}

 Given a sublocale $S$ of $L$, its \emph{closure}\index{closure of a sublocale}\index{sublocale!closure of a $\sim$} and \emph{interior}\index{interior of a sublocale}\index{sublocale!interior of a $\sim$} are  defined by
$$
\overline{S}=\tbigwedge^{\mathsf{S}(L)} \{\mathfrak{c}(a)\mid\mathfrak{c}(a)\ge S\}=\mathfrak{c}(\tbigwedge S)
\quad\text{
and}\quad
S^\circ=\tbigvee^{\mathsf{S}(L)} \{\mathfrak{o}(a)\mid S\ge
\mathfrak{o}(a)\}.
$$

\subsubsection{Subfitness in terms of sublocales}A frame is subfit if and only if each open sublocale is a join of closed sublocales. This was in fact how subfitness was originally defined in \cite{JRI72}.

\subsubsection{Subspaces and induced sublocales} Each subspace $A$ of a $T_0$-space $X$ induces a sublocale of $\mathcal{O}X$ that comprises elements $U\in\mathcal{O}X$ of the form
\[
U=\tbigcup\{V\in\mathcal{O}X\mid V\cap A=U\cap A\}.
\]
We will denote this sublocale by $\mathfrak{s}(A)$. In particular given an open set $U$ one has
\[
\mathfrak{s}(U)=\mathfrak{o}(U)\quad\text{and}\quad
\mathfrak{s}(X\setminus U)=\mathfrak{c}(U).
\]
In general, given $A_i\subseteq X$ for $i\in I$, one has
\[
\mathfrak{s}(\tbigcup_{i\in I}A_i)=\tbigvee_{i\in I}^{\mathsf{S}(L)}\mathfrak{s}(A_i).
\]
Consequently, suprema of induced sublocales are also induced.  Infima of induced sublocales, in contrast, are not necessarily induced. Indeed, we do not necessarily have $\mathfrak{s}(A\cap B)=\mathfrak{s}(A)\wedge\mathfrak{s}(B)$.
 
\subsubsection{Booleanization}\label{booleanization} We can associate with each frame $L$ a complete Boolean algebra $\mathfrak{B}(L)$ consisting of all elements $a=a^{\ast\ast}$ and a frame homomorphism $\beta\colon L\to\mathfrak{B}(L)$ that maps each element $a$ to its double pseudocomplement $a^{\ast\ast}$. $\mathfrak{B}(L)$ is a sublocale of $L$. A sublocale $S$ is said to be \emph{dense} if $\overline{S}=1$, equivalently if $0\in S$. Isbell's density theorem states that each frame $L$ has a least dense sublocale, namely $\mathfrak{B}(L)$.

\subsubsection{The frame $\protect\coS(L)$}
We will make the system of all sublocales of a locale $L$ into a frame $\coS(L)=\mathsf{S}(L)^{\mathrm{op}}$ by considering the dual ordering: $S_{1}\leq S_{2}$ iff $S_{2}\subseteq S_{1}$. Thus, $\{1\}$ is the top and $L$ is the bottom in $\coS(L)$ that we simply denote by $1$ and $0$, respectively. 

\begin{remarkn} In what follows, we will be interested mainly in the frame $\coS(L)$ as it will play an important role in the modeling of discontinuous functions. In consequence, we will consider sublocales as elements of this formal frame almost always and, unless specifically stated otherwise, the order and the lattice operations considered will be always those from $\coS(L)$. Accordingly, for example, we will omit superscripts in joins  and meets, except in ambiguous situations. Also note also that we will be interested in dense elements $S$ of the frame $\coS(L)$, which should not be confused with dense sublocales as in described \ref{booleanization}.
\end{remarkn}

 \subsection{Continuous (extended) real functions}\label{subsection2.2}
Using the basic homomorphism ${\varrho}\colon {\mathfrak{L}\big(\overline{{{\mathbb{R}}}}\big)}\to {{\mathfrak{L}({{\mathbb{R}}})}}$ determined on generators by
 \begin{equation}\label{assignid}
(r,\gidoia)\mapsto(r,\gidoia)\quad\text{and}\quad(\gidoia,s)\mapsto(\gidoia,s)
 \end{equation}
 for each $r,s\in\Q$, the $f\in {{\mathrm{C}(L)}}$ are in a one-to-one correspondence with the $g\in\overline{\mathrm{C}}(L)$ that turn defining relations (r5) and (r6) into identities in $L$ (just take $g=f\varrho$). In what follows we will keep the notation $\mathrm{C}(L)$ to denote also the class inside $\overline{\mathrm{C}}(L)$ of the $f$'s that preserve (r5) and (r6).

$\mathrm{C}(L)$ and $\overline{\mathrm{C}}(L)$ are partially ordered by
\begin{equation}\label{fleg}
\begin{aligned}
f\le g& \iff f(r,\gidoia)\le g(r,\gidoia)\quad\mbox{for all }r\in{{\mathbb{Q}}}\\
&\iff g(\gidoia,s)\le f(\gidoia,s)\quad\mbox{for all }s\in{{\mathbb{Q}}}.
\end{aligned}
\end{equation}

 \begin{examplesn}
For each $r\in \mathbb{Q}$, we have the \emph{constant} function $\boldsymbol{r}\in \mathrm{C}(L)$, given by
\[\boldsymbol{r}(p,{\gidoia})=\begin{cases}0&\text{if }p\ge r
\\ 1&\text{if }p< r \end{cases}
\quad\text{and}\quad
\boldsymbol{r}({\gidoia},q)=\begin{cases}1&\text{if }q>r
\\ 0&\text{if }q\le r\end{cases}\]
for all $p,q\in\Q$. One can similarly has two extended constant functions $\boldsymbol{+\infty}$ and $\boldsymbol{-\infty}$ which are defined for each $p,q\in \mathbb{Q}$ by
 $$\boldsymbol{+\infty}(p,{\gidoia})=1=\boldsymbol{-\infty}({\gidoia},q)\quad\text{and}\quad\boldsymbol{+\infty}({\gidoia},q)=0=\boldsymbol{-\infty}(p,{\gidoia}),$$
and they are precisely the
top and bottom elements of $\overline{\mathrm{C}}(L)$.
\end{examplesn}

An $f\in\mathrm{C}(L)$ is said to be \emph{bounded}\index{bounded!continuous real function}\index{function!bounded $\sim$} if there exist $p,q\in\mathbb{Q}$ such that $\boldsymbol{p}\le f\le \boldsymbol{q}$. Equivalently, $f$ is said to be bounded if and only if there is some rational $r$ such that $f\left((\gidoia,-r)\vee(r,\gidoia)\right)=0$, that is, $f(-r,r)=1$. We shall denote by ${\mathrm{C}}^*(L)$ the set of all bounded members of $\mathrm{C}(L)$. Obviously, all constant functions are in $\mathrm{C}^\ast(L)$.

 As it is well known, $\mathrm{C}(L)$ and $\mathrm{C}^\ast(L)$ are lattices although in general not Dedekind complete \cite{BH03}. In $\mathrm{C}(L)$ (also in $\mathrm{C}^\ast(L)$, since it is sublattice of $\mathrm{C}(L)$) binary joins are given by
 \[
 (f\vee g)(r,\gidoia)=f(r,\gidoia)\vee g(r,\gidoia)\quad\text{and}\quad (f\vee g)(\gidoia,s)=f(\gidoia,s)\wedge g(\gidoia,s)
 \]
for all $r,s\in\Q$ and binary meets are given by 
\[
(f\wedge g)(r,\gidoia)=f(r,\gidoia)\wedge g(r,\gidoia)\quad\text{and}\quad (f\wedge g)(\gidoia,s)=f(\gidoia,s)\vee g(\gidoia,s)
\]
for all $r,s\in\Q$ (see \cite{GP11} for more details).

%
%
%
%

 \subsection{The frame of (extended) partial reals}
Dropping the relation (r2) from the definition of the frame of reals  yields the frame $\mathfrak{L}(\intervalR)$ of \emph{partial} real numbers. It was introduced in \cite{MGP14} as a pointfree counterpart of the partial real line (also interval-domain) which was proposed by Dana Scott in \cite{S72} as a domain-theoretic model for the real numbers. It is a successful idea that has inspired a number of computational models for real numbers. Frame homomorphisms $\mathfrak{L}(\intervalR)\to L$ are called \emph{continuous partial} real functions \cite{MGP14} on $L$. Similarly, dropping (r2) from the definition of the frame of extended reals yields the frame $\mathfrak{L}(\overline{\intervalR})$ of \emph{extended partial} real numbers \cite{MThesis}. Frame homomorphisms $\mathfrak{L}(\overline{\intervalR})\to L$ are called \emph{continuous extended partial} real functions \cite{MThesis} on $L$.

The sets
\[
\mathrm{IC}(L)\quad\text{and}\quad\overline{\mathrm{IC}}(L)
\]
 of continuous partial real functions on $L$ and continuous extended partial real functions on $L$, respectively, are partially ordered by $f\le g$ iff
 \begin{equation}\label{2.5.1}
 f(r,\gidoia)\le g(r,\gidoia)\quad\mbox{and}\quad g(\gidoia,s)\le f(\gidoia,s)
 \end{equation}
  for every $r,s\in\Q$. 
  
  We will also say that an $f\in\mathrm{IC}(L)$ is \emph{bounded} if there exist $p,q\in\Q$ such that $\boldsymbol{p}\leq f\leq\boldsymbol{q}$ and we will denote by $\mathrm{IC}^\ast(L)$ the class of all bounded functions of $\mathrm{IC}(L)$.
  
  \begin{remarks}\label{r2.6} 
 (1) The functions $ h\in \mathrm{IC}(L)$ that factor through the canonical homomorphism ${\iota}\colon \mathfrak{L(\mathbb{IR})}\to\mathfrak{L}(\R)$, determined by the assignment (\ref{assignid}), are just those that turn the defining relation (r2) into an identity in $L$, that is, those which satisfy $ h(r,\gidoia)\vee  h(\gidoia,s)=1$ whenever $r<s$. In view of this, we will keep the notation $\mathrm{C}(L)$ to denote also the class inside $\mathrm{IC}(L)$ of the functions $h$ such that $h(r,\gidoia)\vee h(\gidoia,s)=1$ whenever $r<s$.
 
 The assignment (\ref{assignid}) also determines the canonical homomorphisms  $\overline{\varrho}\colon\mathfrak{L}(\overline{\mathbb{IR}})\to\mathfrak{L}(\mathbb{IR})$ and $\overline{\iota}\colon\mathfrak{L}(\overline{\mathbb{IR}})\to\mathfrak{L}(\overline{\R})$.
An analogous argument motivates us to keep  $\mathrm{IC}(L)$  to denote the class inside $\overline{\mathrm{IC}}(L)$ of the functions $h$ turn the defining relations (r5) and (r6) into identities in $L$ and $\overline{\mathrm{C}}(L)$ to denote the class inside $\overline{\mathrm{IC}}(L)$ of the functions $h$ turn the defining relation  (r2) into identities in $L$.
\smallskip

\noindent (2) In case $f\in\overline{\mathrm{C}}(L)$, as in (\ref{fleg}), the second condition on $f$ and $g$ in (\ref{2.5.1}) is needless because it follows from the first one:
\[
\begin{aligned}g(\gidoia,r)&=g(\tbigvee_{s<r}(\gidoia,s))=\tbigvee_{s<r}g(\gidoia,s)\\
&\le \tbigvee_{s<r}g(s,\gidoia)^*\le
\tbigvee_{s<r}f(s,\gidoia)^*\le f(\gidoia,r),\end{aligned}
\]
 the last inequality because $f$ being in $\overline{\mathrm{C}}(L)$ then, by (r2), $f(s,\gidoia)\vee f(\gidoia,r)=1$ (a similar argument shows that the first condition follows from the second one whenever $g\in \overline{\mathrm{C}}(L)$ and so the two conditions are equivalent if both $f,g$ are in $\overline{\mathrm{C}}(L)$, as in (\ref{fleg})).
 \smallskip

 \noindent (3) There is a dual isomorphism ${-}(\cdot )\colon\overline{\mathrm{IC}}(L)\rightarrow \overline{\mathrm{IC}}(L)$ defined by
\[
({-}h)(\gidoia,r)=h({-}r,\gidoia)\quad\text{ and }\quad ({-}h)(r,\gidoia)=h(\gidoia,{-}r)\quad\text{ for all }r\in
\mathbb{Q}.
\]
When restricted to $A=\mathrm{C}(L),\overline{\mathrm{C}}(L)$ or $\mathrm{IC}(L)$ it yields a dual isomorphism $A\to A$.
 \end{remarks}
 
 \begin{examples} \label{excharfunctions} For each $a,b\in L$ such that $a\wedge b=0$ let $\chi_{a,b}$ denote the bounded continuous partial real function given by
\[
\chi_{a,b}(r,\gidoia)=
 \begin{cases}0&\text{if }r\ge 1,\\a&\text{if }0\le r< 1,\\ 1&\text{if }r< 0, \end{cases}
\quad\text{ and }\quad\chi_{a,b}(\gidoia,s)=
 \begin{cases}1&\text{if }s>1,\\ b&\text{if }0<s\le1,\\ 0&\text{if }s\le 0, \end{cases}
\]
for each $r,s\in{{\mathbb{Q}}}$. Similarly, let $\overline{\chi}_{a,b}$ denote the 
continuous extended partial real function given by
\[
\overline{\chi}_{a,b}(r,\gidoia)=a\quad\text{and}\quad\overline{\chi}_{a,b}(\gidoia,s)=b
\]
for each $r,s\in\Q$. Clearly, $\chi_{a,b}\in\mathrm{C}^\ast(L)$ and $\overline{\chi}_{a,b}\in\overline{\mathrm{C}}(L)$ if and only if $a\vee b=1$, i.e. if and only if $a$ is
complemented with complement $b$.
 \end{examples}

\section{Semicontinuous functions on frames}\label{sectionsemicontinuous}

A lower resp. upper semicontinuous function on a space $X$ is a continuous map $X\to\R_u$ resp. $X\to \R_l$, where $\R_u$ resp. $\R_l$ denotes the space of real numbers with the upper topology resp. lower topology. One obtain their extended versions by replacing the real numbers by the extended real numbers, that is, an extended lower resp. extended upper semicontinuous on $X$ is a continuous map $X\to \overline{\R}_u$ resp. $X\to \overline{\R}_l$, where $\overline{\R}_u$ resp. $\overline{\R}_l$ denotes the space of extended real numbers with the upper topology resp. lower topology. We will denote by $\mathrm{LSC}(X)$, $\overline{\mathrm{LSC}}(X)$, $\mathrm{USC}(X)$ and $\overline{\mathrm{USC}}(X)$ the classes of lower semicontinuous functions on $X$, extended lower semicontinuous functions on $X$, upper semicontinuous functions on $X$ and extended semicontinuous functions on $X$, respectively.

In this section we will first analyze the pointfree counterpart of extended semicontinuous functions and then introduce a conservative definition of semicontinuous functions. 

 \subsection{Extended semicontinuous functions}
 Let $\mathfrak{L}_u(\overline{\R})$ be the subframe of $\mathfrak{L}(\overline{\R})$ generated by elements $(r,\gidoia)$. The frame of open subsets $\overline{\R}_u$ is isomorphic to $\mathfrak{L}_u(\R)$. Since $\overline{\R}_u$ is sober, one also has that $\Pt\mathfrak{L}_u(\overline{\R})\simeq \overline{\R}_u$. By the adjoint situation between frames and topological spaces we have the natural isomorphisms
 \[
\mathsf{Frm}(L,\O X)\simeq \mathsf{Top}(X,\Pt L).
 \]
For $L=\mathfrak{L}_u(\overline{\R})$ one obtains
\[
\mathsf{Frm}(\mathfrak{L}_u(\overline{\R}), \O X)\simeq \mathsf{Top}(X,\overline{\R}_u)=\overline{\mathrm{LSC}}(X).
\]
Accordingly, regarding frame homomorphisms $\mathfrak{L}_u(\overline{\R})\to L$ as extended lower semicontinuous functions on a general frame $L$ provide a conservative extension of the classical notion. One can argue dually to show that regarding frame homomorphisms $\mathfrak{L}_l(\overline{\R})\to L$, where $\mathfrak{L}_l(\overline{\R})$ is the subframe of $\mathfrak{L}(\overline{\R})$ generated by elements $(\gidoia,s)$, as extended upper semicontinuous functions on $L$ extends the classical notion.

The $f$ in $\mathsf{Frm}(\mathfrak{L}_u(\overline{\R}), L)$ resp. $\mathsf{Frm}(\mathfrak{L}_l(\overline{\R}), L)$ are in one-to-one correspondence with the $g\in\overline{\mathrm{C}}(\coS(L))$ such that $g(r,\gidoia)\in\mathfrak{c}L$ for all $r\in\Q$ resp. $g(\gidoia,s)\in\mathfrak{c}L$ for all $s\in\Q$ \cite{GKP2009}.
This motivates the following definition. We will say that $g\in\overline{\mathrm{C}}(\coS(L))$ is \emph{extended lower semicontinuous} functions on $L$ if $g(r,\gidoia)$ is closed for all $r\in\Q$ and we will say that $g$ is \emph{extended upper semicontinuous} functions on $L$ if $g(\gidoia,s)$ is closed for all $s\in\Q$. We will denote by $\overline{\mathrm{LSC}}(L)$ the set of all extended lower semicontinuous functions on $L$ and by $\overline{\mathrm{USC}}(L)$ the set of all extended upper semicontinuous functions on $L$. Furthermore, the $f$ in $\overline{\mathrm{C}}(L)$ are in one-to-one correspondence with the $g$ in $\overline{\mathrm{C}}(\coS(L))$ such that $f(r,\gidoia),f(\gidoia,s)\in\mathfrak{c}L$ for all $r,s\in\Q$. Accordingly, we will keep $\overline{\mathrm{C}}(L)$ to denote the class $\overline{\mathrm{LSC}}(L)\cap\overline{\mathrm{USC}}(L)$.
 
\subsection{Semicontinuous functions}
More complicated is to find  proper pointfree counterparts of the non-extended versions of this notions. Indeed, as pointed out in \cite{GP07}, since the space $\R_u$ is not sober, lower  semicontinuous functions defined on $X$ are not properly represented by frame homomorphism $\mathfrak{L}_u(\R)\to L$, where $\mathfrak{L}_u(\R)$  is the subframe of $\LR$ generated by elements $(r,\gidoia)$. One has that $\Pt\mathfrak{L}_u(\R)$ is homeomorphic to the topological space $\R_{+\infty}$ with set of points $\R\cup\{+\infty\}$ endowed with the upper topology. Consequently, one has an isomorphism
\[
\Omega\colon\mathsf{Top}(X,\R_{+\infty})\to\mathsf{Frm}(\mathfrak{L}_u(\R),\O X)
\]
where $\Omega(\varphi)\colon\mathfrak{L}_u(\R)\to\O(\R)$ is determined by
\[
\Omega(\varphi)(r,\gidoia)=\varphi^{-1}(r,+\infty]
\]
for each $r\in\Q$. Given $f\in\mathsf{Frm}(\mathfrak{L}_u(\R),\O X)$, one has
\[
\Omega^{-1}(f)(x)=\tbigvee\{r\in\Q\mid x\in f(r,\gidoia)\}
\]
for each $x\in X$.

In \cite[Corollary 4.3]{GP07} (see also \cite{GKP2009}) the authors claimed that $\mathrm{LSC}(X)$ is isomorphic, via the restriction of $\Omega$, with those frame homomorphisms $g\colon\mathfrak{L}_u(\R)\to\O X$ such that
\begin{equation}\label{conditionJavi}
\tbigvee_{r\in\Q}^{\coS(\mathcal{O}X)}\mathfrak{o}(g(r,\gidoia))=1.
\end{equation}
(Note that in \cite{GP07} open frame congruences $\Delta_{g(r,\gidoia)}$ are used instead of sublocales.) However, as the following examples shows, this condition is too restrictive for a faithful representation of all real lower semicontinuous functions:
\begin{example}\label{counterexLSC}
Let $\varphi\colon\Q\to\R_{+\infty}$ be a one-to-one map such that $\varphi(\Q)\subseteq\N$. Note that $\varphi$ is lower semicontinuous, as, for all $r\in\Q$, $\varphi^{-1}(r,+\infty]$ is cofinite and, consequently, open. Further, as $\varphi^{-1}(r,+\infty)$ is also dense, one has that so is $\mathfrak{o}(\varphi^{-1}(r,+\infty])=\mathfrak{s}(\varphi^{-1}(r,\gidoia,+\infty])$ as a sublocale. In consequence, by Isbell's density theorem, 
\[
\tbigvee_{r\in\Q}\mathfrak{o}(\Omega(\varphi)(r,\gidoia))=\tbigvee_{r\in\Q}\mathfrak{o}(\varphi^{-1}(r,+\infty])\neq 1
\]
since, for all $r\in\Q$, one has $\mathfrak{o}(\varphi^{-1}(r,+\infty])\leq\mathfrak{B}(\O\Q)\neq 1$ in $\coS(\R_{+\infty})$.
\end{example}

Condition (\ref{conditionJavi}) seems to be an attempt to reflect the fact that $\varphi\in\mathsf{Top}(X,\R_{+\infty})$ takes values in $\R$ if and only if
\[
\tbigcap_{r\in\Q}\varphi^{-1}(r,+\infty]=\varnothing.
\]
As the previous example shows, this is not the case. After all,  in $\coS(L)$, joins of induced sublocales are not necessarily induced. On the other hand, again in $\coS(L)$, a meet of sublocales induced by subspaces is the sublocale induced by the union of those subspace. This suggests an alternative approach, as $\varphi\in\mathsf{Top}(X,\R_{+\infty})$ takes values in $\R$ if and only if
\[
\tbigcup_{r\in \Q}\varphi^{-1}(-\infty,r]=X.
\]

\begin{proposition}\label{isosemicontinuous} Let $X$ be a $T_D$-space.
The restriction of $\Omega$ yields an isomorphism between $\mathrm{LSC}(X)$ and 
\[
\mathcal{A}=\left\{f\colon\mathfrak{L}_u(\R)\to \O X\in\mathsf{Frm}\mid \tbigwedge_{r\in\Q}^{\coS(\mathcal{O}X)}\mathfrak{c}(f(r,\gidoia))=0\right\}.
\]
\end{proposition}

\begin{proof}
For each $\varphi\in\mathrm{LSC}(X)$, one has
\[
\begin{aligned}
\tbigwedge_{r\in\Q}\mathfrak{c}(\Omega(\varphi)(r,\gidoia))&=\tbigwedge_{r\in\Q}\mathfrak{c}(\varphi^{-1}(r,+\infty])\\
&=\tbigwedge_{r\in\Q}\mathfrak{s}(\varphi^{-1}(-\infty, r])\\
&=\mathfrak{s}(\tbigcup_{r\in\Q}\varphi^{-1}(-\infty,r])\\
&=\mathfrak{s}(X)=\mathcal{O}X=0_{\coS(\O X)}.
\end{aligned}
\]

On the other hand,  let $f\colon\mathfrak{L}_u(\R)\to\O X\in\mathsf{Frm}$ such that $\tbigwedge_{r\in\Q}\mathfrak{c}(f(r,\gidoia))=0.$
Since
\[
\begin{aligned}
\mathfrak{s}\left(\tbigcup_{r\in\Q}X\setminus f(r,\gidoia)\right)&=\tbigwedge_{r\in\Q}\mathfrak{s}(X\setminus f(r,\gidoia))\\
&=\tbigwedge_{r\in\Q}\mathfrak{c}(f(r,\gidoia))\\
&=0_{\coS(\O X)}=\O X,
\end{aligned}
\]
one has
\[
\tbigcup_{r\in\Q}X\setminus f(r,\gidoia)=X.
\]
Consequently, for each $x\in X$  there exists some $r\in\Q$ such that $x\not\in f(r,\gidoia)$. Hence, $\Omega^{-1}(f)(x)\leq r$ by (r5). Accordingly, $\Omega(f)(x)\in\R$ for all $x\in X$.
\end{proof}

\begin{remark}
As considered in \cite{GKP2009}, having $\coS(L)$ as common codomain is a convenient approach that allows us to consider lower and upper semicontinuous functions in a common setting. However, note that while
\[
\tbigvee_{r\in\Q}\mathfrak{o}(g(r,\gidoia))=1
\]
implies
\[
\tbigwedge_{r\in\Q}\mathfrak{c}(g(r,\gidoia))=0,
\]
as
\[
\left(\tbigvee_{r\in\Q}\mathfrak{o}(g(r,\gidoia))\right)^\ast=\tbigwedge_{r\in\Q}\mathfrak{c}(g(r,\gidoia)),
\]
Example \ref{counterexLSC} shows that the converse implication does not hold in general. In consequence, this general setting where both lower and upper semicontinuous functions can be defined has to include extended real functions $\mathfrak{L}(\overline{\R})\to\coS(L)$. By the isomorphism $L\simeq \mathfrak{c}(L)^{\mathrm{op}}$, the following definitions properly generalize the classical notion of semicontinuous functions.
\end{remark}

\subsection{}\label{defsemicontinuous}
Given a frame $L$, we will say that
\begin{enumerate} \item a \emph{lower semicontinuous function} on a $L$ is a frame homomorphism $f\colon\mathfrak{L}(\overline{\R})\to \coS(L)$ such that
\begin{enumerate}[(l1)] \item$f(r,\gidoia)\in \mathfrak{c}L$ for all $r\in \Q$, 
\item$\tbigvee_{r\in\Q}f(r,\gidoia)=1$
\item $\tbigwedge_{r\in\Q}f(r,\gidoia)=0$.
\end{enumerate}
\smallskip
\item an \emph{upper semicontinuous function} on $L$ is a frame homomorphism $f\colon\mathfrak{L}(\overline{\R})\to \coS(L)$ such that
\begin{enumerate}[(u1)] \item$f(\gidoia,s)\in \mathfrak{c}L$ for all $s\in \Q$, 
\item$\tbigvee_{s\in\Q}f(\gidoia,s)=1$
\item $\tbigwedge_{s\in\Q}f(\gidoia,s)=0$.
\end{enumerate}
\end{enumerate}
\smallskip
We will denote by $\mathrm{LSC}(L)$ the class of all lower semicontinuous functions on $L$ and by $\mathrm{USC}(L)$ the class of all upper semicontinuous functions on $L$.

As in the case of continuous extended real functions, by the isomorphism between $L$ and $\mathfrak{c}(L)^{\mathrm{op}}$, the $f\in\mathrm{C}(L)$ are in one-to-one correspondence with the $g\in\mathrm{C}(\coS(L))$ such that $g(r,\gidoia),g(\gidoia,s)\in\mathfrak{c}(L)$ for all $r,s\in\Q$. Accordingly, in what follows we will keep $\mathrm{C}(L)$  to denote the class of $g$ in $\mathrm{C}(\coS(L))$ such that $g(r,\gidoia)$ and $g(\gidoia,s)$ are closed for all $r,s\in\Q$. Note that 
\[
\mathrm{C}(L)=\mathrm{LSC}(L)\cap\mathrm{USC}(L).
\]
\begin{remark}\label{dualisosemi}
It is straighforward to check that the restriction of the dual isomorphism in \ref{r2.6}(3) yields dual isomorphims $\mathrm{LSC}(L)\to\mathrm{USC}(L)$ and $\overline{\mathrm{LSC}}(L)\to\overline{\mathrm{USC}}(L)$.\end{remark}

\section{The Dedekind-MacNeille completion of $\overline{\mathrm{C}}(L)$}

In \cite{GMP14,GMP15a,MGP14} several representations of the Dedekind completions of the lattices $\mathrm{C}(L)$ and $\mathrm{C}^\ast(L)$ (see also \cite{MThesis}) were presented. For the aim of this paper we will need a description of Dedekind-MacNeille completion of the lattice of continuous extended real functions $\overline{\mathrm{C}}(L)$ and we will derive it from the Dedekind completion of $\mathrm{C}^\ast(L)$ presented in \cite{MGP14}.

\begin{proposition}\label{proposition4.0.2}
$\overline{\mathrm{IC}}(L)$ is isomorphic to 
\[
\{g\in\mathrm{IC}(L)\mid \boldsymbol{-1}\leq g\leq\boldsymbol{1}\}.
\]
\end{proposition}

\begin{proof}
Let $\alpha\colon\Q\to \Q\cap (-1,1)$ be an order isomorphism with inverse $\beta$. Given $f\in\overline{\mathrm{IC}}(L)$. Let $\Psi(f)\colon\mathfrak{L}(\mathbb{IR})\to L$ be a frame homomorphism determined on  generators by
\[
\Psi(f)(r,\gidoia)=\begin{cases}
0&\text{if }1\leq r\\
f(\beta(r),\gidoia)&\text{if }-1<r<1\\
\tbigvee_{t\in\Q}f(t,\gidoia)&\text{if }r=-1\\
1&\text{if }r<-1
\end{cases}
\]
and
\[
\Psi(f)(\gidoia,s)=\begin{cases}
1&\text{if }1<s\\
\tbigvee_{t\in\Q}f(\gidoia,t)&\text{if }s=1\\
f(\gidoia,\beta(s))&\text{if }-1< s<1\\
0&\text{if }s\leq-1
\end{cases}
\]
for each $r,s\in\Q$. One can easily check that this assignment turns the defining relations (r1), (r3)--(r6) into identities in $L$, consequently, it actually determines a frame homomorphism. Obviously, $\boldsymbol{-1}\leq \Psi(f)\leq \boldsymbol{1}$. Therefore, $f\mapsto\Psi(f)$ defines a map
\[
\Psi\colon\overline{\mathrm{IC}}(L)\to\{g\in\mathrm{IC}(L)\mid \boldsymbol{-1}\leq g\leq\boldsymbol{1}\}.
\]
One can easily check that $\Psi$ is monotone.

Dually, given $g\in\mathrm{IC}(L)$ such that $\boldsymbol{-1}\leq g\leq \boldsymbol{1}$, let $\Phi(g)\colon\mathfrak{L}(\overline{\mathbb{IR}})\to L$ be a frame homomorphism determined on generators by
\[
\Phi(g)(r,\gidoia)=
g(\alpha(r),\gidoia)
\quad
\text{and}
\quad
\Phi(g)(\gidoia,s)=
g(\gidoia,\alpha(s))
\]
for each $r,s\in\Q$. Obviously, these assignments turn the defining relation (r1), (r3) and (r4) into identities in $L$. Therefore, $f\mapsto\Phi(g)$ defines a map
\[
\Phi\colon\{g\in\mathrm{IC}(L)\mid \boldsymbol{-1}\leq g\leq\boldsymbol{1}\}\to\overline{\mathrm{IC}}(L)
\]
which is obvioulsy monotone. Finally, we shall check that $\Psi$ and $\Phi$ are inverse to each other. Let $f\in \overline{\mathrm{IC}}(L)$. Then one has
\[
\Phi\Psi(f)(r,\gidoia)=\Psi(f)(\alpha(r),\gidoia)=f(\beta\alpha(r),\gidoia))= g(r,\gidoia)
\]
for each $r\in\Q$. Dually, one can check that $\Phi\Psi(f)(\gidoia, s)=f(\gidoia,s)$ for each $s\in\Q$. We conclude that $\Phi\Psi(f)=f$. On the other hand, let $g\in\mathrm{IC}(L)$ such that $\boldsymbol{-1}\leq g\leq \boldsymbol{1}$. Then for each $r\in\Q$ such that $1< r$ or $r<-1$ one trivially has $\Psi\Phi(g)(r,\gidoia)=f(r,\gidoia)$. If $-1< r<1$, one has
\[
\Psi\Phi(g)(r,\gidoia)=\Phi(\beta(r),\gidoia)=g(\alpha\beta(r),\gidoia)=g(r,\gidoia),
\]
Finally, for $r=-1$, one has
\[
\begin{aligned}
\Psi\Phi(g)(-1,\gidoia)&=\tbigvee_{t\in\Q}\Phi(g)(t,\gidoia)\\
&=\tbigvee_{t\in\Q}g(\alpha(t),\gidoia)\\
&=\tbigvee_{-1<r<1}g(t,\gidoia)\\
&=g(-1,\gidoia)
\end{aligned}
\]
by (r3). Analogously, one can check that $\Psi\Phi(g)(\gidoia,s)=g(\gidoia,s)$ for all $s\in\Q$. Consequently, $\Phi^{-1}=\Psi$.
\end{proof}

\subsection{Hausdorff continuous functions}
We will say that $f\in\overline{\mathrm{IC}}(L)$ is \emph{Hausdorff continuous} if
\[
f(r,\gidoia)^\ast\leq f(\gidoia,s)\text{ and }f(\gidoia,s)^\ast\leq f(r,\gidoia)\text{ for all } r<s \text{ in }\Q. 
\]
We will denote by $\overline{\mathrm{H}}(L)$ the family of all Hausdorff continuous functions  in $\overline{\mathrm{IC}}(L)$ and by $\mathrm{H}(L)$ the family $\mathrm{IC}(L)\cap\overline{\mathrm{H}}(L)$. It was shown in \cite{MGP14} that, for a complete regular frame $L$,
\[
\mathrm{H}^\ast(L)=\mathrm{IC}^\ast(L)\cap\mathrm{H}(L)
\]
is a Dedekind complete lattice and the inclusion of $\mathrm{C}^\ast(L)$ into $\mathrm{H}^\ast(L)$ is its Dedekind completion.

\begin{remarks}\label{remarksisos}
(1) In \cite{MGP14}, the family $\mathrm{H}^\ast(L)$ was denoted by $\mathrm{C}^\ast(L)^{\vu}$. The motivation for our notation comes from the fact that functions in $\mathrm{H}^\ast(L)$ are the pointfree counterpart of bounded \emph{Hausdorff continuous functions} (see \cite{RA04,Danet2011}).
\smallskip

\noindent (2) Given $f\in\mathrm{IC}(L)$, one has that $f\in{\mathrm{H}}^\ast(L)$ if and only if $\Phi(f)\in\overline{\mathrm{H}}(L)$. In fact, one has that
\[
\Phi(f)(r,\gidoia)^\ast=f(\alpha(r),\gidoia)^\ast\leq f(\gidoia, \alpha(s))=\Phi(f)(\gidoia,s)
\]
for all $r<s$ in $\Q$ if and only if $g(p,\gidoia)^\ast\leq g(q,\gidoia)$ for all $p<q$ in $\Q$. One can check the other condition dually. Consequently, $\overline{\mathrm{H}}(L)$ is isomorphic to
\[
\{f\in\mathrm{H}^\ast(L)\mid \boldsymbol{-1}\leq f\leq\boldsymbol{1}\}.
\]

\noindent (3) Similarly, $f\in \mathrm{C}(L)$ if and only if $\Phi(f)\in\overline{\mathrm{C}}(L)$. Simply note that one has
\[
\Phi(f)(r,\gidoia)\vee\Phi(f)(\gidoia,s)=f(\alpha(r),\gidoia)\vee f(\gidoia,\alpha(s))=1
\]
for all $r<s$ in $\Q$ if and only if $f(p,\gidoia)\vee f(\gidoia,q)=1$ for all $p<q$ in $\Q$. Consequently, $\overline{\mathrm{C}}(L)$ is isomorphic to
\[
\{g\in\mathrm{C}^\ast(L)\mid \boldsymbol{-1}\leq f\leq\boldsymbol{1}\}.
\]

\noindent (4) It is straightforward to check that the dual isomorphism $-(\cdot)\colon\overline{\mathrm{IC}}(L)\to\overline{\mathrm{IC}}(L)$ in \ref{r2.6} (3) restricted to $\overline{\mathrm{H}}(L)$ yields a dual isomorphism $\overline{\mathrm{H}}(L)\to\overline{\mathrm{H}}(L)$.
\smallskip

\noindent (5) Recall $\chi_{a,b}$ and $\overline{\chi}_{a,b}$ from \ref{excharfunctions}. One has that $\chi_{a,b}$ and $\overline{\chi}_{a,b}$ are Hausdorff continuous if and only if $a=b^\ast$ and $b=a^\ast$, that is, if and only if $a=a^{\ast\ast}$ and $b=a^\ast$.
\end{remarks}

\begin{lemma}\label{formulahausdorff}
For each $f\in\overline{\mathrm{H}}(L)$ and $r,s\in\Q$, one has
\begin{enumerate}[\rm(1)]
\item $f(r,\gidoia)=\tbigvee_{p>r}f(\gidoia,p)^\ast=\tbigvee_{p>r}f(p,\gidoia)^{\ast\ast}$ and
\item $f(\gidoia,s)=\tbigvee_{q<s}f(q,\gidoia)^\ast=\tbigvee_{q<s}f(\gidoia,q)^{\ast\ast}$.
\end{enumerate}
\end{lemma}

\begin{proof}
In order to check (1), simply note that for each $r>p$ in $\Q$, one has
\[
f(p,\gidoia)\leq f(p,\gidoia)^{\ast\ast}\leq f(\gidoia, p)^\ast\leq f(r,\gidoia),
\]
since $f(\gidoia,p)\leq f(p,\gidoia)^\ast$ by (r1). Consequently, by (r3), 
\[
f(r,\gidoia)=\tbigvee_{p>r}f(p,\gidoia)=\tbigvee_{p>r}f(\gidoia,p)^\ast=\tbigvee_{p>r}f(p,\gidoia)^{\ast\ast}.
\]
One can check (2) dually.
\end{proof}

\begin{remark}\label{r4.1.3}
By \ref{formulahausdorff}, the argument in \ref{r2.6} is also valid if $f\in\overline{\mathrm{H}}(L)$ and consequently the second condition on $f$ and $g$ in (\ref{2.5.1}) is needless because it follows from the first one. Dually, also the first condition follows from the second one whenever $g\in \overline{\mathrm{H}}(L)$.

\end{remark}

\subsection{The Dedekind-MacNeille completion of $\overline{\mathrm{C}}(L)$}

\begin{proposition}\label{extendedcompletion}
Let $L$ be a completely regular frame. Then $\overline{\mathrm{H}}(L)$ is the Dedekind-MacNeille completion of $\overline{\mathrm{C}}(L)$.
\end{proposition}

\begin{proof}
Just note that since $\mathrm{H}^\ast(L)$ is the Dedekind completion of $\mathrm{C}^\ast(L)$ then
\[
\{f\in\mathrm{H}^\ast(L)\mid \boldsymbol{-1}\leq f\leq\boldsymbol{1}\}
\]
is the Dedekind completion of 
\[
\{f\in\mathrm{C}^\ast(L)\mid \boldsymbol{-1}\leq f\leq\boldsymbol{1}\}.
\]
By the remark above, we conclude that $\overline{\mathrm{H}}(L)$ is the Dedekind completion of $\overline{\mathrm{C}}(L)$. As $\overline{\mathrm{C}}(L)$ is bounded, it is its Dedekind-MacNeille completion. 
\end{proof}

\begin{remark}\label{remarksisos1}
Given a completely regular frame $L$, one can obtain joins and meets in $\overline{\mathrm{H}}(L)$ from those of $\mathrm{H}^\ast(L)$ (see \cite{MGP14}). Let $\{f_i\}_{i\in I}\subseteq\overline{\mathrm{H}}(L)$. Then for $\tbigvee^{\overline{\mathrm{H}}(L)}_{i\in I}=f_\vee$ one has
\[
f_\vee(r,\gidoia)=\tbigvee_{p>r}\left(\tbigvee_{i\in I}f_i(p,\gidoia)\right)^{\ast\ast}
\quad\text{and}\quad
f_\vee(\gidoia,s)=\tbigvee_{q<s}\left(\tbigvee_{i\in I}f_i(q,\gidoia)\right)^{\ast}
\]
Dually, for $\tbigwedge^{\overline{\mathrm{H}}(L)}_{i\in I}=f_\wedge$, one has
\[
f_\wedge(r,\gidoia)=\tbigvee_{p>r}\left(\tbigvee_{i\in I}f_i(\gidoia,q)\right)^{\ast}
\quad\text{and}\quad
f_\wedge(\gidoia,s)=\tbigvee_{q<s}\left(\tbigvee_{i\in I}f_i(\gidoia,q)\right)^{\ast\ast}.
\]
\end{remark}

\begin{proposition}\label{extHausdorffBoolean}
$\overline{\mathrm{H}}(L)$ is isomorphic to $\overline{\mathrm{C}}(\mathfrak{B}(L))$.
\end{proposition}

\begin{proof}
For each $f\in\overline{\mathrm{H}}(L)$ define $\Gamma(f)=\beta_L\cdot f\colon\mathfrak{L}(\overline{\mathrm{IR}})\to\mathfrak{B}(L)$. Note that $f(r,\gidoia)\vee f(\gidoia,s)$ is dense in $L$ whenever $r<s$, as
\[
(f(r,\gidoia)\vee f(\gidoia,s))^\ast=f(r,\gidoia)^\ast\wedge f(\gidoia,s)^\ast\leq f(r,\gidoia)^\ast\wedge f(r,\gidoia)=0.
\]
Consequently,
\[
f(r,\gidoia)^{\ast\ast}\vee_{\mathfrak{B}(L)}f(\gidoia,s)^{\ast\ast}\geq\left(f(r,\gidoia)\vee_L f(\gidoia,s)\right)^{\ast\ast}=1.
\]
Thus $\Gamma(f)\in \overline{\mathrm{C}}(\mathfrak{B}(L))$. Obviously, the map $\Gamma\colon\overline{\mathrm{H}}(L)\to\overline{\mathrm{C}}(\mathfrak{B}(L))$ is order-preserving.

On the other hand, for each $g\in\overline{\mathrm{C}}(\mathfrak{B}(L))$, let $\Delta(g)\colon\mathfrak{L}(\overline{\mathbb{IR}})\to L$ be the frame homomorphism determined on generators as follows:
\[
\Delta(g)(r,\gidoia)=\tbigvee^L_{p>r}g(p,\gidoia)\quad\text{and}\quad \Delta(g)(\gidoia,s)=\tbigvee^L_{q<s}g(\gidoia,q).
\]
It is straightforward to check that those assignments turn the defining relations (r1), (r3) and (r4) into identities in $L$. Let $r<s$ in $\Q$. For each $p,q\in\Q$ such that $r<p<q<s$, one has
\[
\Delta(g)(r,\gidoia)^\ast\leq g(p,\gidoia)^\ast\leq g(\gidoia,q)\leq\Delta(g)(\gidoia,s).
\]
Similarly, $\Delta(g)(\gidoia,s)^\ast\leq\Delta(g)(r,\gidoia)$. Therefore one has a map  $\Delta\colon\overline{\mathrm{C}}(\mathfrak{B}(L))\to\overline{\mathrm{H}}(L)$. It is straightforward to check that it is order-preserving.

Finally, for each $f\in\overline{\mathrm{H}}(L)$, by Lemma \ref{formulahausdorff}, one has
\[
\Delta(\Gamma(f))(r,\gidoia)=\tbigvee^L_{p>r}\Gamma(f)(p,\gidoia)=\tbigvee^L_{p>r}f(p,\gidoia)^{\ast\ast}=f(r,\gidoia)
\]
for each $r\in \Q$. Dually, one can check that $\Delta(\Gamma(f))(\gidoia,s)=f(\gidoia,s)$ for all $s\in\Q$. For each $g\in\overline{\mathrm{C}}(\mathfrak{B}(L))$, one has
\[
\begin{aligned}
\Gamma(\Delta(g))(r,\gidoia)&=\Delta(g)(r,\gidoia)^{\ast\ast}\\
&=\left(\tbigvee^L_{p>r}g(r,\gidoia)\right)^{\ast\ast}\\
&=\tbigvee_{p>r}^{\mathfrak{B}(L)}g(r,\gidoia)=g(r,\gidoia).
\end{aligned}
\]
Thus $\Delta\cdot\Gamma=1_{\overline{\mathrm{H}}(L)}$ and $\Gamma\cdot\Delta=1_{\overline{\mathrm{C}}(\mathfrak{B}(L))}$.
\end{proof}

\begin{corollary}
For each $L$ completely regular frame, $\overline{\mathrm{C}}(\mathfrak{B}(L))$ is isomorphic to the Dedekind-MacNeille completion of $\overline{\mathrm{C}}(L)$.
\end{corollary}

Finally, we close this section by a corollary that extends \cite[Corollary 4.8]{MGP14}.

\begin{proposition}\label{eqextdisc} For any completely regular frame $L$, the following are equivalent:
\begin{enumerate}[\rm(1)]
\item $L$ is extremally disconnected.
\item $\overline{\mathrm{C}}(L)=\overline{\mathrm{H}}(L)$.
\item $\overline{\mathrm{C}}(L)$ is complete.
\end{enumerate}
\end{proposition}

\begin{proof}
(1)$\implies$(2): Let $f\in\overline{\mathrm{H}}(L)$. If $L$ is extremally disconnected, for each $r<t<s$ in $\Q$, one has, by \ref{formulahausdorff},
\[
\begin{aligned}
f(r,\gidoia)\vee f(\gidoia,s)&=\tbigvee_{p>r}f(p,\gidoia)^{\ast\ast}\vee\tbigvee_{q<s}f(s,\gidoia)^\ast\\
&\geq f(t,\gidoia)^{\ast\ast}\vee f(t,\gidoia)^\ast=1.
\end{aligned}
\]
Consequently, $f\in\overline{\mathrm{C}}(L)$.
\smallskip

\noindent (2)$\implies$(1): If $\overline{\mathrm{H}}(L)=\overline{\mathrm{C}}(L)$,   then  $\overline{\chi}_{a^\ast,a^{\ast\ast}}$ is a continuous extended real function on $L$ for each $a\in L$. Consequently,
\[
1=\overline{\chi}_{a^\ast,a^{\ast\ast}}(0,\gidoia)\vee\overline{\chi}_{a^\ast,a^{\ast\ast}}(\gidoia,1)=a^\ast\vee a^{\ast\ast}.
\]
Therefore $L$ is extremally disconnected.
\smallskip

\noindent (2)$\iff$(3): This follows trivially from \ref{extendedcompletion}.
\end{proof}

\section{Arbitrary extended real functions}\label{sectionarbitrary}
The motivation for our approach to arbitrary real functions is based on the following fact: for a $T_1$ space $X$, \emph{$\overline{\mathrm{F}}(X)$ is the smallest complete lattice containing all extended upper and lower semicontinuous functions}, in other words, $\overline{\mathrm{F}}(X)$ the is Dedekind-MacNeille completion of $\overline{\mathrm{LSC}}(X)\cup\overline{\mathrm{USC}}(X)$. We begin by presenting the definition of arbitrary extended real functions in this section.

\subsection{}\label{s5.1}A space $X$ is $T_1$ if and only if $\overline{\mathrm{LSC}}(X)$ is meet-dense in $\overline{F}(X)$ if and only if  $\overline{\mathrm{USC}}(X)$ is join-dense in $\overline{\mathrm{F}}(X)$. Indeed, given $X$ is a $T_1$-space, $x\in X$ and $p\in\R$, the map $g_{x,p}\colon X\to\overline{\R}$ defined as follows is lower semicontinuous:
\[
g_{x,p}(y)=
\begin{cases}
p&\text{if }y=x\\
+\infty&\text{else}
\end{cases}
\]
for each $y\in X$. Given an arbitrary extended real function $f\colon X\to \overline{\R}$, one has that
\[
f=\tbigwedge_{x\in X}g_{x,g(x)}.
\]
Accordingly, $\overline{\mathrm{LSC}}(X)$ is meet-dense in $\overline{\mathrm{F}}(X)$.

On the other hand, if $\overline{\mathrm{LSC}}(X)$ is meet dense in $\overline{\mathrm{F}}(X)$, in particular, given $x\in X$ and $p\in\R$, $g_{x,p}$ is a meet of extended lower semicontinuous functions. In consequence, there exists some $g\in \overline{\mathrm{LSC}}(X)$ such that $g\neq\boldsymbol{+\infty}$ and $g_{x,q}\leq g$. Then $g^{-1}(q,+\infty)=X\setminus\{x\}\in\mathcal{O}X$. Thus $X$ is $T_1$.

By the dual order-isomorphism $-(\cdot)\colon\overline{\mathrm{F}}(X)\to\overline{\mathrm{F}}(X)$ and the fact that it restricts to a dual order isomorphism between $\overline{\mathrm{LSC}}(X)$ and $\overline{\mathrm{USC}}(X)$, we conclude that $X$ is $T_1$ if and only if $\overline{\mathrm{USC}}(X)$ join dense in $\overline{F}(X)$.

\subsection{The  Dedekind-MacNeille completion of $\overline{\mathrm{LSC}}(L)\cup\overline{\mathrm{USC}}(L)$}

For each $a\in L$ and $q\in \Q$, let us denote by $l_{a,q}$ the frame homomorphism $\mathfrak{L}(\R)\to\coS(L)$ determined on generators by the assignment
\[
(r,\gidoia)\mapsto\begin{cases}
\mathfrak{c}(a)&\text{if }r\geq q\\
1&\text{if }r<q
\end{cases}
\quad\text{and}\quad
(\gidoia,s)\mapsto\begin{cases}
\mathfrak{o}(a)&\text{if }s> q\\
0&\text{if }s\leq q
\end{cases}
\]
for each $r,s\in\Q$. Obviously $l_{a,q}\in\overline{\mathrm{LSC}}(L)$.

%

\begin{lemma}
A frame $L$ is subfit if and only if each $S\in\mathfrak{B}(\coS(L))$ is a meet of closed sublocales.
\end{lemma}

\begin{proof}
First let $L$ be a subfit frame. Recall that each $S\in\coS(L)$ can be represented as
\[
\tbigvee_{i\in I}\mathfrak{c}(a_i)\wedge\mathfrak{o}(b_i)
\]
in $\coS(L)$. Since $\mathfrak{o}(a_i)\vee\mathfrak{c}(b_i)$ is complemented for each $i\in I$, one has
\[
\begin{aligned}
S^\ast&=\left(\tbigvee_{i\in I}(\mathfrak{o}(a_i)\vee\mathfrak{c}(b_i))^\ast\right)^\ast\\
&=\tbigwedge_{i\in I}(\mathfrak{o}(a_i)\vee\mathfrak{c}(b_i))^{\ast\ast}\\
&=\tbigwedge_{i\in I}\mathfrak{o}(a_i)\vee\mathfrak{c}(b_i).
\end{aligned}
\]
If $L$ is subfit, for each $i\in I$, there exists $\{d_j\}_{j\in J_i}\subseteq L$ such that
\[
\mathfrak{o}(a_i)=\tbigwedge_{j\in J_i}\mathfrak{c}(d_j).
\]
Since, for each $i\in I$, $\mathfrak{c}(b_i)$ is complemented, one has
\[
\left(\tbigwedge_{J\in J_i}\mathfrak{c}(d_j)\right)\vee\mathfrak{c}(b_i)=\tbigwedge_{j\in J_i}(\mathfrak{c}(d_j)\vee\mathfrak{c}(b_i))=\tbigwedge_{j\in J_i}\mathfrak{c}(d_j\vee b_i )
\]
Therefore
\[
S^\ast=\tbigwedge_{i\in I}\tbigwedge_{j\in J_i}\mathfrak{c}(d_j\vee b_i).
\]
In consequence, all sublocales in $\mathfrak{B}(\coS(L))$ are meets of closed sublocales.

On the other hand, if each $S\in\coS(L)$ is a meet of closed sublocales, in particular, so are all open sublocales as they are complemented.
\end{proof}

\begin{proposition}
Let $L$ be a subfit frame. Then
$\overline{\mathrm{LSC}}(L)$ is meet-dense in $\overline{\mathrm{C}}(\coS(L))$ and $\overline{\mathrm{USC}}(L)$ is join-dense in $\overline{\mathrm{C}}(\coS(L))$.
\end{proposition}

\begin{proof}
Let $f\in\overline{\mathrm{C}}(\coS(L))$ and
\[
A_q=\{a\in L\mid\mathfrak{c}(a)\geq f(q,\gidoia)^{\ast\ast}\}.
\]
for each $q\in\Q$. Note that, as $L$ is subfit, one has
\[
f(q,\gidoia)^{\ast\ast}=\tbigwedge_{a\in A_q}\mathfrak{c}(a).
\]
We shall show that
\[
f=\tbigwedge^{\overline{\mathrm{H}}({\scriptsize\coS}(L))}_{q\in\Q,a\in A_q}l_{a,q}.
\]
Recall that $\overline{\mathrm{H}}(\coS(L))$ is the Dedekind-MacNeille completion of $\overline{\mathrm{C}}(\coS(L))$, therefore complete. By \ref{remarksisos1}, for each $r\in\Q$, one has
\[
\begin{aligned}
\left(\tbigwedge^{\overline{\mathrm{H}}({\scriptsize\coS}(L))}_{q\in\Q,a\in A_q}l_{a,q}\right)(r,\gidoia)&=\tbigvee_{s>r}\left(\tbigvee_{q\in\Q}\tbigvee_{a\in A_q}l_{a,q}(\gidoia,s)\right)^\ast\\
&=\tbigvee_{s>r}\tbigwedge_{q\in\Q}\tbigwedge_{a\in A_q}l_{a,q}(\gidoia,s)^\ast\\
&=\tbigvee_{s>r}\tbigwedge_{q<s}\tbigwedge_{a\in A_q}\mathfrak{c}(a).
\end{aligned}
\]
Since $A_q\subseteq A_s$ whenever $q<s$, we conclude that
\[
\left(\tbigwedge^{\overline{\mathrm{H}}({\scriptsize \coS}(L))}_{q\in\Q,a\in A_q}l_{a,q}\right)(r,\gidoia)=\tbigvee_{s>r}\tbigwedge_{a\in A_s}\mathfrak{c}(a)=\tbigvee_{s>r}f(s,\gidoia)^{\ast\ast}=f(r,\gidoia)
\]
and accordingly
\[
f=\tbigwedge^{\overline{\mathrm{H}}({\scriptsize \coS}(L))}_{q\in\Q,a\in A_q}l_{a,q}.
\]
Therefore $\overline{\mathrm{LSC}}(L)$ in meet-dense in $\overline{\mathrm{C}}(\coS(L))$. The fact that $\overline{\mathrm{USC}}(L)$ is join-dense in $\overline{\mathrm{C}}(\coS(L))$ follows easily from the dual isomorphism $-(\cdot)\colon\overline{\mathrm{C}}(\coS(L))\to \overline{\mathrm{C}}(\coS(L))$.
\end{proof}

\begin{corollary}\label{semicontinuouscompletion}
Let $L$ be a subfit frame. Then the Dedekind-MacNeille completions of $\overline{\mathrm{LSC}}(L)\cup\overline{\mathrm{USC}}(L)$ and $\overline{\mathrm{C}}(\coS(L))$ coincide. 
\end{corollary}

\subsection{Arbitrary extended real function on $L$}  \ref{extendedcompletion} and \ref{semicontinuouscompletion}, combined with the fact that $\coS(L)$ is always completely regular, motivate to define \emph{arbitrary extended real functions on a frame $L$} as extended Hausdorff continuous functions on $\coS(L)$, that is, an arbitrary extended real function on $L$ is a frame homomorphism
\[
f\colon\mathfrak{L}(\overline{\mathbb{IR}})\to\coS(L)
\]
such that $f(r,\gidoia)^\ast\leq f(\gidoia,s)$ and $f(\gidoia,s)^\ast\leq f(r,\gidoia)$ for all $r<s$. We will denote by
\[
\overline{\mathrm{F}}(L)
\]
the family of all arbitrary extended real functions on a subfit frame $L$

\begin{remark}By \ref{extHausdorffBoolean}, one has that
\[
\overline{\mathrm{F}}(L)\simeq \overline{\mathrm{C}}(\mathfrak{B}(\coS(L))).
\]
\end{remark}

\subsection{Joins and meets in $\overline{\mathrm{F}}(L)$} For the sake of completeness, we provide formulae for joins and meets in $\overline{\mathrm{F}}(L)$, even though this is just a particular case of the lattice operations on the lattice of extended Hausdorff continuous functions on an arbitrary frame described in \ref{remarksisos1}. Let $\{f_i\}_{i\in I}\subseteq\overline{\mathrm{F}}(L)$. Then for $\tbigvee^{\overline{\mathrm{F}}(L)}_{i\in I}=f_\vee$ one has
\[
f_\vee(r,\gidoia)=\tbigvee_{p>r}\left(\tbigvee_{i\in I}f_i(p,\gidoia)\right)^{\ast\ast}
\quad\text{and}\quad
f_\vee(\gidoia,s)=\tbigvee_{q<s}\left(\tbigvee_{i\in I}f_i(q,\gidoia)\right)^{\ast}
\]
for each $r\in \Q$.  Dually, for $\tbigwedge^{\overline{\mathrm{F}}(L)}_{i\in I}=f_\wedge$, one has
\[
f_\wedge(r,\gidoia)=\tbigvee_{p>r}\left(\tbigvee_{i\in I}f_i(\gidoia,p)\right)^{\ast}
\quad\text{and}\quad
f_\wedge(\gidoia,s)=\tbigvee_{q<s}\left(\tbigvee_{i\in I}f_i(\gidoia,q)\right)^{\ast\ast}.
\]
for each $s\in\Q$.

\section{Arbitrary real functions}
In this section we will address the main goal of the paper: defining arbitrary real functions on a frame $L$. One might be tempted to simply replace $\mathfrak{L}(\overline{\mathbb{IR}})$ by $\mathfrak{L}(\mathbb{IR})$ in the definition of $\overline{\mathrm{F}}(L)$. However, one would want $\mathrm{F}(L)$ to contain all semicontinuous functions but, for instance, a lower semicontinuous function $f\in\mathrm{LSC}(L)$ does not necessarily turn the defining relation (r6) into the identity in $\coS(L)$.

\subsection{Identifying real-valued functions}

An arbitrary extended real function $f\colon X\to \overline{\R}$ is real-valued, that is, $f\in\mathrm{F}(X)$, if and only if 
\begin{enumerate}
\item $f\vee_{\overline{\mathrm{F}}(X)} g=\boldsymbol{+\infty}$ 
 implies $g=\boldsymbol{+\infty}$,
\item $f\wedge_{\overline{\mathrm{F}}(X)} g=\boldsymbol{-\infty}$ 
 implies $h=\boldsymbol{-\infty}$.
\end{enumerate}
Accordingly, we will say that $f\in\overline{\mathrm{F}}(L)$ is \emph{real-valued} or an \emph{arbitrary real function} on $L$ if
\begin{enumerate}
\item $f\vee_{\overline{\mathrm{F}}(L)} g=\boldsymbol{+\infty}$ 
  implies $g=\boldsymbol{+\infty}$,
\item $f\wedge_{\overline{\mathrm{F}}(L)} g=\boldsymbol{-\infty}$ 
 implies $h=\boldsymbol{-\infty}$.
\end{enumerate}
We will denote by $\mathrm{F}(L)$ the class of all arbitrary real functions on $L$.

\begin{remark} It is straightforward to check that the dual isomorphism $-(\cdot)\colon\overline{\mathrm{F}}(L)\to\overline{\mathrm{F}}(L)$ from \ref{remarksisos} (4) when restricted to $\mathrm{F}(L)$ yields a dual isomorphism $\mathrm{F}(L)\to\mathrm{F}(L)$.
\end{remark}

\begin{proposition}
Let $L$ be a frame and $f\in\overline{\mathrm{F}}(L)$. The following are equivalent:
\begin{enumerate}[\rm(1)]
\item $f\vee^{\overline{\mathrm{F}}(L)} g=\boldsymbol{+\infty}$ 
 implies $g=\boldsymbol{+\infty}$,
\item $\left(\tbigvee_{s\in\Q}f(\gidoia,s)\right)^\ast=0$,
\item $\tbigwedge_{s\in\Q}f(r,\gidoia)=0$.
\end{enumerate}

\end{proposition}

\begin{proof}
(1)$\implies$(2): Let $f\in\overline{\mathrm{F}}(L)$ such that $f\vee^{\overline{\mathrm{F}}(L)} g=\boldsymbol{+\infty}$  implies $g=\boldsymbol{+\infty}$. Let $S\in\coS(L)$. Recall $\overline{\chi}_{S^{\ast\ast},S^\ast}$ from \ref{excharfunctions}. Obviously, $\overline{\chi}_{S^{\ast\ast},S^\ast}\in\overline{\mathrm{F}}(L)$ and $\overline{\chi}_{S^{\ast\ast},S^\ast}=+\infty$ if and only if $S^{\ast\ast}=1$. Consequently, if $S^{\ast\ast}\neq 1$, one has $f\vee^{\overline{\mathrm{F}}(L)}\overline{\chi}_{S^{\ast\ast},S^\ast}\neq +\infty$. Then one has
\[
\begin{aligned}
0&\neq \tbigvee_{s\in\Q}\left(f\vee^{\overline{\mathrm{F}}(L)}\overline{\chi}_{S^{\ast\ast},S^\ast}\right)(\gidoia,s)\\
&=\tbigvee_{s\in\Q}\tbigvee_{q<s}(f(q,\gidoia)\vee\overline{\chi}_{S^{\ast\ast},S^\ast}(q,\gidoia))^\ast\\
&=\tbigvee_{s\in\Q}f(s,\gidoia)^\ast\wedge\overline{\chi}_{S^{\ast\ast},S^\ast}(s,\gidoia)^\ast\\
&=\tbigvee_{s\in\Q}f(s,\gidoia)^\ast\wedge S^\ast\\
&=S^\ast\wedge\tbigvee_{s\in\Q}f(s,\gidoia)^\ast.
\end{aligned}
\]
Therefore, by \ref{formulahausdorff},
\[
S^\ast\wedge \tbigvee_{s\in\Q}f(\gidoia,s)\neq 0
\]
whenever $S^{\ast\ast}\neq 1$, equivalently, whenever $S^\ast\neq 0$. We conclude that
\[
\left(\tbigvee_{s\in\Q}f(\gidoia,s)\right)^\ast=0.
\]

\noindent (2)$\iff$(3): This follows easily from the fact that 
\[
 f(\gidoia, t)^\ast\leq f(s,\gidoia)\leq f(\gidoia,s)^\ast
\]
for all $s<t$ in $\Q$ and
\[
\left(\tbigvee_{s\in\Q}f(\gidoia,s)\right)^\ast=\tbigwedge_{s\in\Q}f(\gidoia, s)^\ast.
\]

\noindent (2)$\implies$(1): Let $g\in \overline{\mathrm{F}}(L)$ such that $f\vee g=\boldsymbol{+\infty}$. For each $q\in\Q$, one has
\[
0=(f\vee g)(\gidoia, q)=\tbigvee_{s<q}(f(s,\gidoia)\vee g(s,\gidoia))^\ast=\tbigvee_{s<q}(f(s,\gidoia)^\ast\wedge g(s,\gidoia)^\ast).
\]
Then $g(s,\gidoia)^\ast\leq f(s,\gidoia)^{\ast\ast}$ for each $s\in\Q$, as $f(s,\gidoia)^\ast\wedge g(s,\gidoia)^\ast=0$. Moreover, for each $r<s$ in $\Q$ one has
\[
g(r,\gidoia)^\ast\leq f(s,\gidoia)^{\ast\ast}
\]
since $g(r,\gidoia)\geq g(s,\gidoia)$. One also has 
\[
g(s,\gidoia)^\ast\leq f(r,\gidoia)^{\ast\ast}
\]
as $f(s,\gidoia)\leq f(r,\gidoia)$. In consequence, 
\[
\begin{aligned}
g(s,\gidoia)^\ast&\leq\tbigwedge_{r\in\Q}f(r,\gidoia)^{\ast\ast}\\
&=\left(\tbigvee_{r\in\Q}f(r, \gidoia)^\ast\right)^\ast\\
&\leq\left(\tbigvee_{r\in \Q}f(\gidoia,r)\right)^\ast=0.
\end{aligned}
\]
Then, by (r1), one has
\[
g(\gidoia,s)\leq g(s,\gidoia)^\ast=0
\]
for all $s\in\Q$. We conclude that $g=\boldsymbol{+\infty}$.
\end{proof}

By a dual argument, we conclude the folllowing.
\begin{corollary}\label{charreal}
Let $L$ be a frame and $f\in\overline{\mathrm{F}}(L)$. The following are equivalent:
\begin{enumerate}[\rm(1)]
\item $f$ is real-valued,
\item $\left(\tbigvee_{r\in\Q}f(r,\gidoia)\right)^\ast=0=\left(\tbigvee_{s\in\Q}f(\gidoia,s)\right)^\ast$,
\item $\tbigwedge_{s\in\Q}f(\gidoia,s)=0=\tbigwedge_{r\in\Q}f(r,\gidoia)$.
\end{enumerate}
\end{corollary}

\begin{proposition}
For any frame $L$, 
\[
\mathrm{LSC}(L)=\mathrm{F}(L)\cap\overline{\mathrm{LSC}}(L)\quad\text{and}\quad\mathrm{USC}(L)=\mathrm{F}(L)\cap\overline{\mathrm{USC}}(L).
\]
\end{proposition}

\begin{proof}
Let $f\in\mathrm{F}(L)\cap\overline{\mathrm{LSC}}(L)$.  In order to show that $f\in\mathrm{LSC}(L)$ we only have to check that condition (l2) in \ref{defsemicontinuous} holds. For each $r\in\Q$, there exists $a_r$ such that $f(r,\gidoia)=\mathfrak{c}(a_r)$. Then
\[
0=\left(\tbigvee_{r\in\Q}f(r,\gidoia)\right)^\ast=\tbigwedge_{r\in\Q}\mathfrak{o}(a_r)=\mathfrak{o}(\tbigvee_{r\in\Q}a_r).
\]
Accordingly, $\tbigvee_{r\in\Q}a_r=1$. We conclude that
\[
\tbigvee_{r\in\Q}f(r,\gidoia)=\tbigvee_{r\in\Q}\mathfrak{c}(a_r)=\mathfrak{c}\left(\tbigvee_{r\in\Q}a_r\right)=\mathfrak{c}(1)=1.
\]
Consequently, one has
\[
\mathrm{F}(L)\cap\overline{\mathrm{LSC}}(L)\subseteq\mathrm{LSC}(L).
\]
The reverse inclusion is straightforward. One can check dually that $\mathrm{USC}(L)=\mathrm{F}(L)\cap\overline{\mathrm{USC}}(L)$.
\end{proof}

\subsection{Conservativeness} 
Recall that the poset of extended semicontinuous functions on a space is isomorphic to the poset of extended semicontinuous functions on its frame of open sets. Taking into account that for a $T_1$-space $X$ one has that  $\mathcal{O}X$ is subfit and that the Dedekind-MacNeille completion of a poset is unique up to isomorphism, we obtain
\[
\overline{\mathrm{F}}(X)\simeq\overline{\mathrm{F}}(\mathcal{O}X).
\]
Furthermore, as the definition of real-valued functions only relies on the order structure of $\overline{\mathrm{F}}(\mathcal{O}X)$, one can easily conclude that
\[
\mathrm{F}(X)\simeq\mathrm{F}(\mathcal{O}X).
\]

\subsection{Semicontinuous regularizations}We shall show now how to extend lower and upper regularizations studied in \cite{GKP09a} to our new setting (see also \cite{GKP2009}). The \emph{lower regularization} $f^\circ\in\overline{\mathrm{LSC}}(L)$ of $f\in\overline{\mathrm{F}}(L)$ is determined on generators by
\[
f^\circ(r,\gidoia)=\tbigvee_{p>r}\overline{f(p,\gidoia)}\quad\text{and}\quad f^\circ(\gidoia,s)=\tbigvee_{q<s}\left(\overline{f(q,\gidoia)}\right)^\ast
\]
for each $r,s\in\Q$. It is straightforward to check that this assignment turn the defining relations (r1)--(r4) into identities in $\coS(L)$. Further, if $f$ is real-valued, then
\[
\tbigwedge_{r\in\Q}f^\circ(r,\gidoia)\leq\tbigwedge_{r\in\Q}f(r,\gidoia)=0,
\]
thus (l3) holds. However (l2) does not hold in general: $f^\circ\in\mathrm{LSC}(L)$ if and only if
\[
\tbigvee_{r\in\Q}\overline{f(r,\gidoia)}=1.
\]

Dually, the \emph{upper regularization} $f^-\in\overline{\mathrm{USC}}(L)$ of $f\in\overline{\mathrm{F}}(L)$ is defined by
\[
f^-=-(-f)^\circ.
\]
An easy computation gives
\[
f^-(r,\gidoia)=\tbigvee_{p>r}\left(\overline{f(\gidoia,p)}\right)^\ast\quad\text{and}\quad f^-(\gidoia,s)=\tbigvee_{q<s}\overline{f(\gidoia,q)}
\]
for each $r,s\in\Q$. Further, one has that $f^-\in\mathrm{USC}(L)$ if $f\in\mathrm{F}(L)$ and
\[
\tbigvee_{s\in\Q}\overline{f(\gidoia,s)}=1.
\]

One can follow the arguments in \cite[Propositions 7.3 and 7.4]{GKP2009} to show that in our new setting $(\cdot)^\circ\colon\overline{\mathrm{F}}(L)\to\overline{\mathrm{LSC}}(L)$ resp. $(\cdot)^-\colon\overline{\mathrm{F}}(L)\to\overline{\mathrm{USC}}(L)$ is an interior-like operator resp. a closure-like operator.

\subsection{}
An alternative definition for arbitrary real functions on a frame $L$ is considered  in \cite{PP16}  where, in the particular where $L$ is subfit, arbitrary real functions are homomorphism $g\colon\mathfrak{L}(\R)\to \mathfrak{B}(\mathsf{S}(L))$ (see also \cite{PPT16}).  Even though in \cite{PP16} sublocales are ordered by inclusion, we can equivalently keep the dual order as in the rest of this paper, since $\mathfrak{B}(\mathsf{S}(L))$ is a Boolean algebra, thus consequently it is dually isomorphic to itself. Recall the isomorphism $\Gamma\colon\overline{\mathrm{F}}(L)\to\overline{\mathrm{C}}(\mathfrak{B}(\coS(L)))$ from \ref{extHausdorffBoolean}. Note that, for $f\in\overline{\mathrm{F}}(L)$, one has
\[
1=\tbigvee_{r\in\Q}^{\mathfrak{B}({\scriptsize\coS}(L))}\Gamma(f)(r,\gidoia)=\left(\tbigvee^{{\scriptsize\coS}(L)}_{r\in\Q}f(r,\gidoia)^{\ast\ast}\right)^{\ast\ast}
\]
if and only if
\[
0=\left(\tbigvee^{{\scriptsize\coS}(L)}_{r\in\Q}f(r,\gidoia)^{\ast\ast}\right)^\ast=\tbigwedge^{{\scriptsize\coS}(L)}_{r\in\Q}f(r,\gidoia)^\ast=\left(\tbigvee^{{\scriptsize\coS}(L)}_{r\in\Q}f(r,\gidoia)\right)^\ast.
\]
Analogously, one has
\[
\tbigvee_{s\in\Q}^{\mathfrak{B}({\scriptsize\coS}(L))}\Gamma(f)(\gidoia,s)=1
\]
if  and only if
\[
\left(\tbigvee^{{\scriptsize\coS}(L)}_{s\in\Q}f(\gidoia,s)\right)^\ast=0.
\]
Consequently, the isomorphism $\Gamma$ restrists to an isomorphism
\[
\mathrm{F}(L)\to \mathrm{C}(\mathfrak{B}(\coS(L))),
\]
showing that both approaches to the notion of arbitrary real functions are equivalent for subfit frames. 

Furthermore, this shows that we can define algebraic operations in such a way that $\mathrm{F}(L)$ becomes a lattice-ordered ring. Consequently, there is no impediment to develop a theory of rings of real functions in this setting.

\section{$\mathrm{F}(L)$ vs. $\mathrm{C}(\protect\coS(L))$}
We conclude this article by analyzing the relation between our approach to arbitrary real functions and the one adopted in \cite{GKP2009}. The following diagram summarizes the relations between the classes of functions considered in this paper (each arrow represents an inclusion, which is strict in the general case):
 \begin{center}
 \begin{tikzpicture}
 \matrix (m) [matrix of math nodes, row sep=3em,
 column sep=2em]{
 &  & \overline{\mathrm{F}}(L) &   \\
 & \overline{\mathrm{LSC}}(L) & \overline{\mathrm{C}}(\coS(L)) &   \overline{\mathrm{USC}}(L) \\
 & \mathrm{F}(L)  & \overline{\mathrm{C}}(L) &  \\
\mathrm{LSC}(L)  & \mathrm{C}(\coS(L)) & \mathrm{USC}(L) &  \\
 & \mathrm{C}(L) &  &  \\};
 \path[-stealth]
 (m-5-2) edge (m-4-1)
 (m-5-2) edge (m-4-2)
 (m-5-2) edge (m-4-3)
 (m-5-2) edge (m-3-3)
 (m-4-1) edge (m-3-2)
 (m-4-1) edge (m-2-2)
 (m-4-2) edge (m-3-2)
 (m-4-3) edge (m-2-4)
 (m-3-3) edge (m-2-2)
 (m-3-3) edge (m-2-3)
 (m-3-3) edge (m-2-4)
 (m-2-2) edge (m-2-3)
 (m-2-3) edge (m-1-3)
 (m-2-4) edge (m-2-3)
 (m-4-2) edge[-,line width=6pt,draw=white] (m-2-3) edge (m-2-3)
 (m-4-3) edge[-,line width=6pt,draw=white] (m-3-2) edge (m-3-2)
 (m-3-2) edge[-,line width=6pt,draw=white] (m-1-3) edge (m-1-3);
\end{tikzpicture}
\end{center}
The only inclusion that we have not explicitly considered yet, $\mathrm{C}(\coS(L))\subseteq\mathrm{F}(L)$, follows from \ref{charreal}, since
\[
\tbigvee_{r\in\Q}f(r,\gidoia)=1=\tbigvee_{s\in\Q}f(\gidoia,s)
\]
trivially implies
\[
\left(\tbigvee_{r\in\Q}f(r,\gidoia)\right)^\ast=0=\left(\tbigvee_{s\in\Q}f(\gidoia,s)\right)^\ast.
\]

The difference between the extended and non-extended cases has already been analyzed in \ref{sectionsemicontinuous} and \ref{sectionarbitrary}. By \ref{eqextdisc}, we know that
\[
\overline{\mathrm{F}}(L)=\overline{\mathrm{C}}(\coS(L))
\]
if and only if $\coS(L)$ is extremally disconnected. More complicated is the case of real-valued functions. By the chararcterization of Banaschewski-Hong \cite[Proposition 1]{BH03}, one has that $\mathrm{C}(\coS(L))$ is Dedekind complete if and only if $\coS(L)$ is extremally disconnected. Since $\coS(L)$ is not extremally disconnected in general, it follows from the following proposition that
\[
\mathrm{F}(L)\neq\mathrm{C}(\coS(L))
\]
in general.

\begin{proposition}
For any frame $L$, $\mathrm{F}(L)$ is Dedekind complete.
\end{proposition}

\begin{proof}
Let $\{f_i\}_{i\in I}\subseteq \mathrm{F}(L)$ and $f\in\mathrm{F}(L)$ such that $f_i\leq f$ for all $i\in I$. As $\overline{\mathrm{F}}(L)$ is complete, there exists
\[
f_\vee=\tbigvee_{i\in I}^{\overline{\mathrm{F}}(L)}f_i\in \overline{\mathrm{F}}(L).
\]
As $f_\vee\leq f$, one has $f_\vee(r,\gidoia)\leq f(r,\gidoia)$ for all $r\in\Q$ and consequently
\[
\tbigwedge_{r\in\Q}f_\vee(r,\gidoia)\leq\tbigwedge_{r\in\Q}f(r,\gidoia)=0.
\]
As $f_i\leq f$ for each $i\in I$, one has $f(\gidoia,s)\leq f_i(\gidoia,s)$ for all $s\in\Q$. Therefore
\[
\tbigwedge_{s\in\Q}f_\vee(\gidoia,s)\leq\tbigwedge_{s\in\Q}f_i(\gidoia,s)=0.
\]
We conclude that $f_\vee$ is real-valued. Since $\mathrm{F}(L)$ is dually isomorphic to itself, then $\mathrm{F}(L)$ is Dedekind complete.
\end{proof}

\begin{example}
Let $\mathfrak{s}(\Q)$ resp. $\mathfrak{s}(\mathbb{I})$ be the sublocale of $\mathcal{O}\R$ induce by the subspace $\Q$ of all rational points resp. by the subspace $\mathbb{I}$ of all irrational points. One can check that $\mathfrak{s}(\Q)$ and $\mathfrak{s}(\mathbb{I})$ are pseudocomplement to each other in $\coS(\mathcal{O}\R)$, that is, $\mathfrak{s}(\Q)^\ast=\mathfrak{s}(\mathbb{I})$ and $\mathfrak{s}(\mathbb{I})^\ast=\mathfrak{s}(\Q)$. Therefore, $\chi_{\mathfrak{s}(\Q),\mathfrak{s}(\mathbb{I})}$ is a Hausdorff continuous partial real function, that is, $\chi_{\mathfrak{s}(\Q),\mathfrak{s}(\mathbb{I})}\in\overline{\mathrm{F}}(L)$. Moreover, it is obviously real-valued. However, it is not in $\mathrm{C}(\coS(\mathcal{O}\R))$, as
\[
\chi_{\mathfrak{s}(\Q),\mathfrak{s}(\mathbb{I})}(0,\gidoia)\vee \chi_{\mathfrak{s}(\Q),\mathfrak{s}(\mathbb{I})}(\gidoia,1)=\mathfrak{s}(\Q)\vee\mathfrak{s}(\mathbb{I})\leq \mathfrak{B}(\mathcal{O}\R)\neq 1,
\]
since both $\mathfrak{s}(\Q)$ and $\mathfrak{s}(\mathbb{I})$ are dense sublocales and consequently
\[
\mathcal{O}\R\neq \mathfrak{B}(\mathcal{O}\R)\subseteq
\mathfrak{s}(\Q)\vee\mathfrak{s}(\mathbb{I}).
\]
\end{example}
\smallskip

For the  following proposition, recall that an element $a\in L$ is said to be a \emph{cozero} if there exists $f\in\mathrm{C}(L)$ such that $a=f(\gidoia,0)\vee f(0,\gidoia)$. In that case, there exists $g\in\mathrm{C}(L)$ such that $a=g(\gidoia,1)$ (simply take $g=((-f)\wedge f)+\boldsymbol{1}$). See \cite{BB97} for more details. A frame $L$ is said to be a \emph{P-frame} if each cozero element is complemented and said to be an \emph{almost P-frame} if $a=a^{\ast\ast}$ for each cozero element. Obvioulsy, if $L$ is an almost P-frame, then $1$ is the only dense cozero. In fact, this is a sufficient condition \cite{Dube2009}. Further, under extremal disconnectedness, $L$ is a P-frame iff and only if it is an almost P-frame.

\begin{proposition}
For any frame $L$,  $\mathrm{F}(L)=\mathrm{C}(\coS(L))$ if and only if $\coS(L)$ is an extremally disconnected P-frame.
\end{proposition}

\begin{proof}
$\implies$:
For each $S\in\coS(L)$, one has $\chi_{S^\ast,S^{\ast\ast}}\in\mathrm{F}(L)$. If $\mathrm{F}(L)=\mathrm{C}(\coS(L))$, one has that
\[
1=\chi_{S^\ast,S^{\ast\ast}}(0,\gidoia)\vee \chi_{S^\ast,S^{\ast\ast}}(\gidoia,1)=S^\ast\vee S^{\ast\ast}.
\]
Consequently, $\coS(L)$ is extremally disconnected.

On the other hand, let $T$ be a dense cozero element of $\coS(L)$. Then there exists $g\in\mathrm{C(\coS(L))}$ such that $T=g(\gidoia,1)$. Le $g^\prime=(g\vee\boldsymbol{0})\wedge\boldsymbol{1}$. Obviously, $\boldsymbol{1}\leq g^\prime\leq \boldsymbol{1}$. By the isomorphism $\Phi$ in \ref{proposition4.0.2}, one has $\Phi(g^\prime)\in\overline{\mathrm{F}}(L)$. Then one has
\[
\begin{aligned}
\left(\tbigvee_{r\in\Q}\Phi(g^\prime)(r,\gidoia)\right)^\ast&=\left(\tbigvee_{-1<r<1} g^\prime(r,\gidoia)\right)^\ast\\
&\leq g^\prime(-1/2,\gidoia)^\ast=1^\ast=0,
\end{aligned}
\]
since $g^\prime\geq\boldsymbol{0}$, and
\[
\begin{aligned}
\left(\tbigvee_{s\in\Q}\Phi(g^\prime)(\gidoia,s)\right)^\ast&=\left(\tbigvee_{-1<s<1} g^\prime(\gidoia,s)\right)^\ast\\
&=g^\prime(\gidoia,1)^\ast\\
&=g(\gidoia,1)^\ast=T^\ast=0,
\end{aligned}
\]
since $T$ is dense. Consequently, $\Phi(g^\prime)\in\mathrm{F}(L)$. If $\mathrm{F}(L)=\mathrm{C}(\coS(L))$, one has
\[
\begin{aligned}
1&=\tbigvee_{s\in\Q}\Phi(g^\prime)(\gidoia,s)\\
&=\tbigvee_{s\in\Q}g^\prime(\gidoia,\alpha(s))\\
&=\tbigvee_{-1<s<1}g^\prime(\gidoia,s)\\
&=g^\prime(\gidoia,s)=g(\gidoia,1)=T.
\end{aligned}
\]

\noindent $\impliedby$:  Let $f\in\mathrm{F}(L)$. We already know that $\overline{\mathrm{F}}(L)=\overline{\mathrm{C}}(\coS(L))$ if $\coS(L)$ is extremally disconnected, consequently $f$ turns (r2) into an identity in $\coS(L)$, thus we only have to check (r5) and (r6). By the isomorphism $\Psi$ from \ref{proposition4.0.2}, one has $\Psi(f)\in\mathrm{C}(\coS(L))$ and
\[
\Psi(f)(-1,\gidoia)=\tbigvee_{r\in\Q}f(r,\gidoia).
\]
Thus $\tbigvee_{r\in\Q}f(r,\gidoia)$ is a cozero, which is dense by \ref{charreal}. In consequence, $\tbigvee_{r\in\Q}f(r,\gidoia)=1$, thus $f$ turns the defining relation (r5) into a identity  in $\coS(L)$. One can check (r6) dually.
\end{proof}

\section*{Acknowledgments} Thanks are given to Andrew Moshier, Joanne Walter-Waylands, Iraide Mardones Garc\'\i a and Mar\'\i a \'Angeles de Prada Vicente for their insightful suggestions and comments.

\end{document}